\documentclass[11pt,a4paper]{article}
\usepackage{amsfonts}
\usepackage{amsmath,amssymb}
\usepackage{mathrsfs}
\usepackage{hyperref}
\usepackage{graphicx}
\usepackage{float}

\newtheorem{thm}{Theorem}[section]

\newtheorem{lem}[thm]{Lemma}
\newtheorem{prop}[thm]{Proposition}
\newtheorem{defn}[thm]{Definition}

\numberwithin{equation}{section}\allowdisplaybreaks

\def\leq{\leqslant}

\def\leq{\leqslant}
\def\geq{\geqslant}

\begin{document}

\title{\large\bf  Local  Well and Ill Posedness for the Modified KdV Equations \\ in Subcritical Modulation Spaces}

\author{\normalsize \bf Mingjuan Chen\footnote{Corresponding author.}
,\  Boling Guo \\
\footnotesize
\it Institute of Applied Physics and Computational Mathematics, Beijing 100088, PR China \\
\footnotesize
\it Emails: mjchenhappy@pku.edu.cn(M.Chen); gbl@iapcm.ac.cn(B.Guo)\ \  \\
}
\date{} \maketitle

\thispagestyle{empty}
\begin{abstract}
 We consider the Cauchy problem of the modified KdV equation (mKdV)
\begin{align}
 u_{t} + u_{xxx} \pm (u^3)_x = 0,\quad
u(0,x)=u_0(x).
\end{align}
Local well-posedness of this problem is obtained in modulation spaces $M^{1/4}_{2,q}(\mathbb{{R}})$ $(2\leq q\leq\infty)$. Moreover, we show that the data-to-solution map fails to be $C^3$ continuous in $M^{s}_{2,q}(\mathbb{{R}})$ when $s<1/4$.  It is well-known that $H^{1/4}$ is a critical Sobolev space of mKdV so that it is well-posedness in $H^s$ for $s\geq 1/4$ and ill-posed (in the sense of uniform continuity) in $H^{s'}$ with $s'<1/4$. Noticing that $M^{1/4}_{2,q} \subset B^{1/q-1/4}_{2,q}$ is a sharp embedding and $H^{-1/4}\subset B^{-1/4}_{2,\infty}$, our results contains all of the subcritical data in $M^{1/4}_{2,q}$, which contains a class of functions in $H^{-1/4}\setminus H^{1/4}$. \\


\end{abstract}

\section{Introduction}

In this paper we study the Cauchy problem of the  modified Korteweg-de Vries (mKdV) equation on the real line $\mathbb{R}$:
\begin{align}
 u_{t} + u_{xxx} \pm (u^3)_x = 0,\quad
u(0,x)=u_0(x),\quad x\in\mathbb{R},
 \label{mKdV}
\end{align}
where $u=u(x,t)\in \mathbb{R}$ with $(x,t)\in \mathbb{R}^{1+1}$.

The scale invariant homogeneous Sobolev space for mKdV is $\dot{H}^{-1/2}$. That is to say, for any solution $u(x,t)$ of \eqref{mKdV} with initial data $u_0(x)$, the scaling function $u_{\lambda}(x,t):=\lambda u(\lambda x, \lambda^3 t)$ is also a solution of \eqref{mKdV} with initial data $u_{0,\lambda}:=\lambda u_0(\lambda x)$, and  satisfies
\begin{align}\label{scale}
 \|u_{0,\lambda}\|_{\dot{H}^{-1/2}}=\|u_0\|_{\dot{H}^{-1/2}}.
\end{align}
On the other hand, $H^{1/4}$ is the critical Sobolev space of mKdV so that it is globally well-posed in $H^s$ for $s\geq 1/4$ and ill-posed in $H^{s'}$ with $s'<1/4$. The ill-posed result is in the sense that the data-to-solution map fails to be uniformly continuous on a fixed ball in $H^{s'}$ with $s'<1/4$. The local well-posed result for $s
\geq1/4$ by using the contraction method  and ill-posed result for the focusing equation ($+$ sign in front of the nonlinearity) were proved by Kenig, Ponce and Vega, see \cite{KPV93} and \cite{KPV01}, respectively. The local well-posed result was extended to a global one for $s>1/4$ due to Colliander, Keel, Staffilani, Takaoka and Tao by using $I$-method, see \cite{CKSTT03}. The global result for $s=1/4$ was obtained by Guo in \cite{Guo09}. In addition, the ill-posed result for the defocusing equation ($-$ sign in front of the nonlinearity) was obtained by Christ, Colliander and Tao \cite{CCT03}.

Therefore, there is $3/4$ derivative gap between ${H}^{-1/2}$ and ${H}^{1/4}$ for the well-posedness result of mKdV. In order to discover the behavior of the solution out of ${H}^{1/4}$, Gr\"unrock brought in the $\widehat{H^{q'}_s}$ spaces for which the norm is defined by
 $$\|u\|_{\widehat{H^{q'}_s}}:=\|\langle\xi\rangle^s\hat{u}\|_{L^q},\ \ \ \ 1/q+1/q'=1,$$
 and he obtained the local well-posedness of \eqref{mKdV} for data $u_0\in \widehat{H^{q'}_s}(\mathbb{R}),  2\leq q<4,  s\geq s(q):=1/2q$ in \cite{G04}. In 2009, Gr\"unrock and Vega broadened the range of $q$ to $2\leq q<\infty$ by using the trilinear estimates in \cite{G09}. From the scaling point, the spaces $\widehat{H^{q'}_s}$ behave like the Sobolev spaces $H^\sigma$, if $s-1/2+1/q=\sigma$. Thus, they can lower the regularity to $-1/2$ by taking $q$ tending to infinity, but there is no result for $q=\infty$. In this paper we consider the initial data in more general modulation spaces $M^s_{2,q}$, $2\leq q\leq \infty$ (Indeed, $\widehat{H^{q'}_s}\subset M^s_{2,q}$).

Modulation space $M^s_{p,q}$ was introduced by Feichtinger \cite{Fei2} in 1983 and equivalently defined in the following way (cf.\cite{WaHaHu09,WaHu07,WaHaHuGu11,WaHud07}):
\begin{align}
\|f\|_{M^s_{p,q}(\mathbb{R})}= \left(\sum_{k\in \mathbb{Z} } \langle k \rangle^{sq}
\|\Box_k f\|^q_{L^p(\mathbb{R})} \right)^{1/q},  \label{mod-space1}
\end{align}
where $\Box_k = \mathscr{F}^{-1} \chi_{[k-1/2, k+1/2]} \mathscr{F}$, $\mathscr{F}$ ($\mathscr{F}^{-1}$) denotes the (inverse) Fourier transform on $ \mathbb{R}$ ,  $\chi_{E}$ denotes the characteristic function on $E$ and  $\langle k \rangle=(1+|k|^2)^{1/2}$. From Plancherel theorem and H\"older's inequality, we know that $\widehat{H^{q'}_s}\subset M^s_{2,q}$ ($2\leq q\leq \infty$). Moreover, combining the sharp inclusions between Besov and modulation spaces, we have (cf. \cite{SuTo07,WaHaHuGu11})
$$
\widehat{H_{1/4}^{q'}} \subset M^{1/4}_{2,q} \subset B^{1/q-1/4}_{2,q}, \ \  2\leq q\leq \infty,
$$
where the inclusions are optimal. Therefore, our result in which the initial data belongs to $M^{1/4}_{2,\infty}$ can be certainly seen as an improvement. Our main theorem is as follows.

\begin{thm}\label{theorem1}
Let $2\leq q \leq\infty$,  $u_0 \in  M^{1/4}_{2,q}$. Then there exists a time $T>0$ such that mKdV \eqref{mKdV} is locally well posed in $C([0,T];  M^{1/4}_{2,q}) \cap X^{1/4}_{q,A}([0,T])$, where $X^{1/4}_{q,A}$ is defined in next section. Moreover, the regularity index $1/4$ in $M^{1/4}_{2,q}$ is optimal. Specifically, if $s<1/4$, the data-to-solution map in $M^{s}_{2,q}(\mathbb{{R}})$ is not $C^3$ continuous at origin.
\end{thm}

Modulation spaces contain a class of initial data out of the critical Sobolev spaces $H^{s_c}$, for which the nonlinear PDE is well-posed for $s>s_c$ and ill-posed for $s<s_c$. Therefore, solving the nonlinear PDE in modulation spaces has absorbed some researchers' attention, see \cite{BO09,BGOR07,CFS12,CWWW18,CN08,CN08b,CN09,Iw10,KaKoIt12,KaKoIt14,Kato14,Wa13,Wang1}. We will use $U^p$ and $V^p$ spaces in our discussion, since the dual relation and other important properties are ideally to deal with the nonlinearity. $U^p$ and $V^p$ spaces are introduced to solving PDEs by Koch and Tataru, see \cite{CHT12,HaHeKo09,KoTa05,KoTa07}. Combining $U^p$, $V^p$  and modulation spaces, Guo, Ren and the second author have considered the cubic and derivative non-linear Schr\"odinger equation, respectively, see \cite{Gu16,GRW16}.

Let us list some notations. Let $c < 1$, $C>1$  denote positive universal constants, which can
be different at different places; $a\lesssim b$ stands for $a\leq C
b$, $a\sim b$ means that $a\lesssim b$ and
$b\lesssim a$; $a\approx b$ means that $|a-b|\leq C$, $a\gg b$ means that $a>b+C$;  We write $a\wedge b =\min(a,b)$, $a\vee b
=\max(a,b)$;  $p'$ is the dual number of $p \in
[1,\infty]$, i.e., $1/p+1/p'=1$.

\section{Function spaces}

\subsection{Definitions}
 In this subsection, we review some function spaces used to obtain the well-posedness theory for non-linear dispersive equations. $U^p$ spaces were first applied by Koch and Tataru \cite{CHT12,KoTa05,KoTa07,KoTa12}, and $V^p$ spaces are due to Wiener \cite{W79}.

Let $\mathcal{Z}$ be the set of finite partitions $-\infty= t_0 <t_1<...< t_{K-1} < t_K =\infty$. In the following, we consider functions taking values in $L^2:=L^2(\mathbb{R}^d;\mathbb{C})$, but in the general $L^2$ may be replaced by an arbitrary Hilbert space or general Banach space.
\begin{defn}
Let $1\leq p <\infty$. For any $\{t_k\}^K_{k=0} \in \mathcal{Z}$ and $\{\phi_k\}^{K-1}_{k=0} \subset L^2$ with $\sum^{K-1}_{k=0} \|\phi_k\|^p_2=1$, $\phi_0=0$. A step function $a: \mathbb{R}\to L^2$ given by
$$
a= \sum^{K}_{k=1} \chi_{[t_{k-1}, t_k)} \phi_{k-1}
$$
is said to be a $U^p$-atom. All of the $U^p$ atoms is denoted by $\mathcal{A}(U^p)$.   The $U^p$ space is
$$
U^p:=\left\{u= \sum^\infty_{j=1} c_j a_j : \ a_j \in \mathcal{A}(U^p), \ \ c_j \in \mathbb{C}, \ \ \sum^\infty_{j=1} |c_j|<\infty  \right\}
$$
for which the norm is given by
$$
\|u\|_{U^p}:= \inf \left\{\sum^\infty_{j=1} |c_j| : \ \ u= \sum^\infty_{j=1} c_j a_j , \ \  \ a_j \in \mathcal{A}(U^p), \ \ c_j \in \mathbb{C} \right\}.
$$
\end{defn}
\begin{defn}
Let $1\leq p <\infty$. We define $V^p$ as the normed space of all functions $v: \mathbb{R} \to L^2$ such that $\lim_{t\to \pm \infty} v(t)$ exist and for which the norm
$$
\|v\|_{V^p} := \sup_{\{t_k\}^K_{k=0} \in \mathcal{Z}} \left( \sum^K_{k=1} \|v(t_k)-v(t_{k-1})\|^p_{L^2}\right)^{1/p}
$$
is finite, where we use the convention that $v(-\infty) = \lim_{t\to \infty} v(t)$ and $v(\infty)=0$ (here $v(\infty)$ and $\lim_{t\to  \infty} v(t)$ are different notations). Likewise, we denote by $V^p_-$ the subspace of all $v\in V^p$ so that $v(-\infty) =0$. Moreover, we define the closed subspace $V^p_{rc}$ $(V^p_{-,rc})$ as all of the right continuous functions in $V^p$ $(V^p_-)$.
\end{defn}

\begin{defn}
We define
$$
U^p_{A} := e^{-\cdot \partial_x^3} U^p, \ \ \|u\|_{U^p_{A}} = \|e^{ t \partial_x^3} u \|_{U^p},
$$
$$
V^p_{A} := e^{-\cdot \partial_x^3} V^p, \ \ \|u\|_{V^p_{A}} = \|e^{ t \partial_x^3} u \|_{V^p},
$$
and similarly for the definition of $V^p_{rc, A}$, $V^p_{-, A}$, $V^p_{-, rc, A}$.
\end{defn}

\begin{defn}
 Besov type Bourgain's spaces $\dot X^{s, b, q}$ are defined by
$$
\|u\|_{\dot X^{s,b,q}} := \left\| \|\chi_{|\tau-\xi^3|\in [2^{j-1}, 2^j)} |\xi|^{s} |\tau-\xi^3|^{b} \widehat{u}(\tau,\xi) \|_{L^2_{\xi,\tau}}  \right\|_{\ell^q_{j\in \mathbb{Z}}}.
$$
\end{defn}

\begin{defn}
The frequency-uniform localized $U^2$-spaces $X^s_{q}(I)$ and $V^2$-spaces $Y^s_{q} (I)$  are defined by
\begin{align}
 \|u\|_{X^s_q(I)} & = \left(\sum_{\lambda \in I \cap \mathbb{Z}} \langle \lambda\rangle^{sq}\|\Box_\lambda u\|^q_{U^2}\right)^{1/q}, \quad  X^s_q:= X^s_q(\mathbb{R}),  \label{mod-space2}\\
 \|v\|_{Y^s_q(I)} & = \left(\sum_{\lambda \in I \cap \mathbb{Z}} \langle \lambda\rangle^{sq}\|\Box_\lambda v\|^q_{V^2}\right)^{1/q}, \quad  Y^s_q:= Y^s_q(\mathbb{R}),  \label{mod-space2a} \\
\|u\|_{X^s_{q, A}}: & = \| e^{t \partial_x^3}  u\|_{X^s_q}, \ \ \ \|v\|_{Y^s_{q, A}}: = \| e^{t \partial_x^3}  v\|_{Y^s_q}.
\end{align}
\end{defn}

\subsection{Known results}
The following known results about $U^p$ and $V^p$ can be found in \cite{Gu16,HaHeKo09,KoTa05,KoTa12}.
\begin{prop} \label{UVprop1}
{\rm (Embedding)} Let $1\leq p <q <\infty$. We have the following results.
\begin{itemize}
 \item[\rm (i)]  $U^p$ and $V^p$, $V^p_{rc}$, $V^p_{-}$, $V^p_{rc, -}$ are Banach spaces.

\item[\rm (ii)] $U^p\subset V^p_{rc, -} \subset U^q \subset L^\infty (\mathbb{R}, L^2)$. Every $u\in U^p$ is right continuous on $t\in \mathbb{R}$.

\item[\rm (iii)] $V^p \subset V^q$,   $V^p_{-} \subset V^q_{-} $,   $V^p_{rc} \subset V^q_{rc} $,  $V^p_{rc, -} \subset V^q_{rc, -} $.

\item[\rm (iv)] $\dot X^{0, 1/2, 1} \subset U^2_{A} \subset V^2_{A} \subset \dot X^{0, 1/2, \infty}$.
\end{itemize}
\end{prop}

Similar to the Schr\"odinger equation, whose dispersive modulation is $|\tau+\xi^2|$, the mKdV equation's dispersive modulation is $|\tau-\xi^3|$.  By the last inclusion of (iv) in Proposition \ref{UVprop1}, we see that

\begin{lem}[\rm Dispersion Modulation Decay]
Suppose that the dispersion modulation $|\tau-\xi^3| \gtrsim \mu$ for a function $u\in L^2_{x,t}$, then we have
\begin{align}
  \|u \|_{L^2_{x,t} }  \lesssim \mu^{-1/2} \|u\|_{V^2_A}. \label{dispersiondecay}
\end{align}
\end{lem}

\begin{prop} \label{UVprop2}
{\rm (Interpolation)} Let $1\leq p <q <\infty$.  There exists a positive constant $\epsilon(p,q)>0$, such that for any $u\in V^p $ and $M>1$,  there exists a decomposition $u=u_1+u_2$ satisfying
\begin{align}
 \frac{1}{M} \|u_1\|_{U^p} + e^{\epsilon M}  \|u_2\|_{U^q} \lesssim \|u\|_{V^p}. \label{interp}
\end{align}
\end{prop}

   Let $I\subset \mathbb{R}$ be an interval with finite length. For the sake of simplicity, we denote
$$
 u_\lambda = \Box_\lambda u, \ \  u_I = \sum_{\lambda\in I \cap \mathbb{Z}} u_\lambda.
$$
\begin{prop} \label{orthogonality}{\rm(orthogonality in $V^2$)} Take an interval $I\subset \mathbb{R}$, then for $u\in V^2$ the following orthogonality holds:
\begin{align}
\|u_I\|_{V^2}\leq \bigg(\sum_{\lambda\in I\cap\mathbb{Z}} \| u_\lambda\|_{V^2}^2\bigg)^{1/2}.
\end{align}
\end{prop}
\begin{prop} \label{UVprop3}
{\rm (Duality)} Let $1\leq p   <\infty$, $1/p+1/p'=1$.  Then $(U^p)^* = V^{p'}$ in the sense that
\begin{align}
T: V^{p'}  \to (U^p)^*; \ \ T(v)=B(\cdot,v), \label{dual}
\end{align}
is an isometric mapping.  The bilinear form $B: U^p\times V^{p'}$ is defined in the following way: For a partition $\mathrm{t}:= \{t_k\}^K_{k=0} \in \mathcal{Z}$, we define
 \begin{align}
B_{\mathrm{t}} (u,v) = \sum^K_{k=1} \langle u(t_{k-1}), \ v(t_k)-v(t_{k-1})\rangle. \label{dual2}
\end{align}
Here $\langle \cdot, \cdot \rangle$ denotes the inner product on $L^2$. For  any $u\in U^p$, $v\in V^{p'}$, there exists a unique number $B(u,v)$ satisfying the following property. For any $\varepsilon>0$, there exists a partition $\mathrm{t}$ such that
$$
|B(u,v)- B_{\mathrm{t}'} (u,v)| <\varepsilon, \ \ \forall \  \mathrm{t}'    \supset \mathrm{t}.
$$
Moreover,
$$
|B(u,v)| \leq \|u\|_{U^p} \|v\|_{V^{p'}}.
$$
In particular, let $u\in V^1_{-}$ be absolutely continuous on compact interval, then for any $v\in V^{p'}$,
$$
 B(u,v) =\int \langle u'(t), v(t)\rangle dt.
$$
\end{prop}

\begin{prop}{\rm \cite{GRW16}} \label{UVprop4} {\rm (Duality)}
 Let $1\leq q   <\infty$.  Then $(X^s_q)^* = Y^{-s}_{q'} $ in the sense that
\begin{align}
T: Y^{-s}_{q'}   \to (X^s_q)^* ; \ \ T(v)=B(\cdot,v), \label{dualprop4}
\end{align}
is an isometric mapping, where the bilinear form $B(\cdot,\cdot)$ is defined in  Proposition \ref{UVprop3}. Moreover, we have
$$
|B(u,v)| \leq   \|u\|_{X^s_q} \|v\|_{Y^{-s}_{q'}}.
$$
\end{prop}

\section{Basic Estimates}

\begin{lem}{\rm \cite{KPV89}(Strichartz Estimates)} \label{strichartz}
Let ($p,q$) satisfy the admissibility condition
\begin{align}\label{admissbility}
  \frac{2}{p}+\frac{1}{q}=\frac{1}{2},\ \ \ \ 4\leq p\leq \infty, \ \ 2\leq q\leq \infty.
\end{align}
Then
\begin{align}\label{strichartza}
\|D_x^{1/p}e^{-t\partial^3_x}\phi\|_{L^p_tL^q_x}\lesssim \|\phi\|_{L^2}.
\end{align}
In particular, for $N\geq1$,
\begin{align}\label{strichartzb}
 \|P_{N}e^{-t\partial^3_x}\phi\|_{L^8_tL^4_x}\lesssim \langle N\rangle^{-1/8} \|\phi\|_{L^2}.
\end{align}
By testing atoms in $U^8_A$ space, we obtain
\begin{align}\label{strichartzc}
 \|P_{N}u\|_{L^8_tL^4_x}\lesssim \langle N\rangle^{-1/8} \|u\|_{U^8_A}.
\end{align}
\end{lem}

\begin{lem}[\rm Bilinear Estimate] \label{V2decay2}
Suppose that $\widehat{u_0}, \widehat{v_0}$ are localized in some compact intervals $I_1,I_2$ with $dist(I_1, I_2)\gtrsim \lambda$, $dist(I_1, -I_2)\gtrsim \mu$. Then,
\begin{align}\label{bilinear1}
 \|e^{-t\partial_x^3}u_0e^{-t\partial_x^3} v_0 \|_{L^2_{x, t}} \lesssim (\lambda\mu)^{-1/2}\|u_0\|_{L^2} \|v_0\|_{L^2}.
\end{align}
By testing atoms in $U^2_A$ space, we obtain
\begin{align}\label{bilinear2}
 \|uv\|_{L^2_{x, t}} \lesssim (\lambda\mu)^{-1/2}\|u\|_{U^2_A} \|v\|_{U^2_A}.
\end{align}
Applying the interpolation in Proposition \ref{UVprop2},  for any $0<\varepsilon \ll 1$ and $0<T\leq 1$, we get
\begin{align}\label{bilinear3}
 \|u v \|_{L^2_{x, t\in[0,T]}} \lesssim T^{\varepsilon/4} (\lambda\mu)^{-1/2 + \varepsilon} \|u\|_{V^2_A} \|v\|_{V^2_A}.
\end{align}
\end{lem}
{\bf Proof.} Taking the Fourier transform in space, we have
\begin{align}
& \mathscr{F}_x  \left( e^{-t\partial_x^3}u_0  e^{-t\partial_x^3}v_0  \right)(\xi,t)  = \int e^{it(\xi^3+3\xi\xi_1^2-3\xi^2\xi_1)} \widehat{u}_0(\xi-\xi_1) \widehat{v}_0(\xi_1) d\xi_1.
 \label{bilinearest1}
\end{align}
Then taking the Fourier transform in time, we obtain
\begin{align}
& \mathscr{F}_{x,t}  \left( e^{-t\partial_x^3}u_0  e^{-t\partial_x^3}v_0  \right)(\xi,\tau)  = \int  \delta(\tau+3\xi^2\xi_1-\xi^3-3\xi\xi_1^2)  \widehat{u}_0(\xi-\xi_1) \widehat{v}_0(\xi_1) d\xi_1.
 \label{bilinearest2}
\end{align}
Denote
$$
g(\xi_1) = \tau+3\xi^2\xi_1-\xi^3-3\xi\xi_1^2,
$$
we see that the zeros and the derivative are
$$
\xi^\pm_1 = \frac{\xi}{2} \pm \sqrt{\frac{\xi^2}{4}- \frac{\xi^3-\tau}{3\xi}}:=\frac{\xi}{2} \pm y, \ \  g'(\xi_1) = 3\xi^2-6\xi\xi_1.
$$
Recall that $\delta(g(\xi_1)) = \delta (\xi_1-\xi_1^+)/ |g'(\xi_1^+)| + \delta (\xi_1-\xi_1^-)/ |g'(\xi_1^-)|= \delta (\xi_1-\xi^+)/ 6|\xi|y +\delta (\xi_1-\xi_1^-)/6|\xi|y $, we have
\begin{align}
 \tiny{\mathscr{F}_{x,t}  \left( e^{-t\partial_x^3}u_0  e^{-t\partial_x^3}v_0  \right)(\xi,\tau)  = \frac{1}{6|\xi|y} \widehat{u}_0 \left(\frac{\xi}{2} - y\right) \widehat{v}_0\left(\frac{\xi}{2} +y\right) +\frac{1}{6|\xi|y} \widehat{u}_0\left(\frac{\xi}{2} + y\right) \widehat{v}_0 \left(\frac{\xi}{2} -y\right).} \label{bilinear4.6}
\end{align}
By symmetry, it suffices to estimate the first term in \eqref{bilinear4.6}.  Changing of variables $y= \sqrt{\frac{\xi^2}{4}- \frac{\xi^3-\tau}{3\xi}}$ and considering $d\tau=c|\xi||y|dy$, we see that
\begin{align}
\left\|  e^{-t\partial_x^3}u_0  e^{-t\partial_x^3}v_0  \right\|^2_{L^2_{x,t}} &  \leq \int_{\mathbb{R}^2}\frac{c}{|\xi||y|}\left |\widehat{u}_0 \left(\frac{\xi}{2} - y \right)\right|^2 \left | \widehat{v}_0\left(\frac{\xi}{2} +y\right)\right|^2 dy d\xi \  \nonumber\\
&  \lesssim \int_{\mathbb{R}^2}\frac{1}{|\xi_1-\xi_2||\xi_1+\xi_2|} |\widehat{u}_0(\xi_1)|^2 | \widehat{v}_0(\xi_2)|^2 d\xi_1d\xi_2  \notag\\
&  \lesssim \lambda^{-1}\mu^{-1} \int_{\mathbb{R}^2} |\widehat{u}_0(\xi_1)|^2 | \widehat{v}_0(\xi_2)|^2 d\xi_1d\xi_2\notag\\
&\lesssim \lambda^{-1}\mu^{-1}\|u_0\|^2_2\|v_0\|^2_2,
\label{bilinearestmm}
\end{align}
where in the last inequality, we have applied $dist(I_1, I_2) \geq\lambda$ and $dist(I_1, -I_2)\geq \mu$. $\hfill\Box$

\begin{lem}\label{L4} {\rm ($L^4$ Estimates)}
Let $I\subset [0, +\infty)$ or $(-\infty, 0]$ with $|I|<\infty$.  For  any  $\theta\in (0,1)$, $\beta>0$, we have
\begin{align}
 \|u_I\|^2_{L^4_{x,t\in [0,T]}} \lesssim (T^{1/4}+ T^{(1-\theta)/4}|I|^{2\beta+(1-\theta)/2 }) \|u\|^2_{X^{-1/8}_{4,A}(I)}. \label{lebesgue4}
\end{align}
In particular, if $1 \lesssim  |I| <\infty$, $0<T<1$, then for any $0< \varepsilon \ll 1$, $4\leq q\leq\infty$
\begin{align}
 \|u_I\|_{L^4_{x,t\in [0,T]}} \lesssim   T^{\varepsilon/4}|I|^{1/4-1/q + \varepsilon} \max_{\lambda \in I} \langle \lambda\rangle^{-3/8} \|u\|_{X^{1/4}_{q,A}(I)}.  \label{lebesgue4a}
\end{align}
\end{lem}
{\bf Proof.} Without loss of generality, we assume $I\subset [0, +\infty)$.
\begin{align*}
\| u_I \|_{L^4([0,T]\times \mathbb{R})}^2 &= \| (u_I)^2 \|_{L^2([0,T]\times \mathbb{R})}=\bigg\| \sum_{m, n\in I\cap \mathbb{Z}} u_m u_n\bigg\|_{L^2([0,T]\times \mathbb{R})} \\
	&\leq \sum_{k\in \mathbb{N}}\bigg\| \sum_{m-n\sim2^k}u_m u_n\bigg\|_{L^2([0,T]\times \mathbb{R})}.
\end{align*}

Case $k=0$, i.e. $m\approx n$:
\begin{align*}
& \bigg\| \sum_{n\in I\cap \mathbb{Z}}u_n^2\bigg\|_{L^2([0,T]\times \mathbb{R})}\leq \bigg(\sum_{n\in I\cap \mathbb{Z}}\| u_n^2\|_{L^2([0,T]\times \mathbb{R})}^2\bigg)^{1/2}\\
	\leq &\  \bigg(\sum_{n\in I\cap \mathbb{Z}} \| u_n\|_{L^4([0,T]\times \mathbb{R})}^4\bigg)^{1/2} \leq T^{1/4} \bigg(\sum_{n\in I\cap \mathbb{Z}} \| u_n\|_{L^8_{t\in[0,T]}L^4_x}^4\bigg)^{1/2}\\
\lesssim &\  T^{1/4} \bigg(\sum_{n\in I\cap\mathbb{Z}}\Big(\langle n\rangle ^{-1/8}\| u_n\|_{U_A^8}\Big)^4\bigg)^{1/2}
\lesssim T^{1/4}\|u\|^2_{X^{-1/8}_{4,A}},
\end{align*}
where the first step is by the orthogonality in $L^2$ and the last step follows from the Strichartz estimate.

 Case $k>0$: Notice that $k$ is summed for $\ln|I|$ times, we have $\sum_{k\in \mathbb{N}}\lesssim \ln|I|$. We split the other sum as follows
$$
\sum_{m-n\sim2^k}u_m u_n=
\sum_{n\in I\cap\mathbb{Z}}\sum\limits_{m\in I\cap \mathbb{Z},\atop m-n\sim 2^k}u_m u_n=
\sum_{j\in \mathbb{N}^+}\sum\limits_{n\in I\cap\mathbb{Z},\atop n\sim j 2^k}\sum\limits_{m\in I\cap\mathbb{Z},\atop m-n\sim 2^k}u_m u_n,
$$
where $j$ is chosen such that $j2^k, (j+1)2^k\in I$. Hence for $u_n$ with $n\sim j 2^k$ and $u_m$ with $m-n\sim 2^k$, we have that the frequency of the function $u_m u_n$ will be close to $(2j+1) 2^k$, which implies by orthogonality that

\begin{align}\label{lebesgue4b}
&\sum_{k\in \mathbb{N}}\bigg \| \sum_{m-n\sim2^k}u_m u_n\bigg\|_{L^2([0,T]\times \mathbb{R})}=\sum_{k\in \mathbb{N}}\bigg \| \sum_{j\in \mathbb{N}^+}\sum\limits_{n\in I\cap\mathbb{Z},\atop n\sim j 2^k}\sum\limits_{m\in I\cap\mathbb{Z},\atop m-n\sim 2^k}u_m u_n\bigg\|_{L^2([0,T]\times \mathbb{R})}\nonumber\\
&\lesssim \sum_{k\in \mathbb{N}}\bigg(\sum_{j\in \mathbb{N}^+} \bigg\|\sum\limits_{n\in I\cap\mathbb{Z},\atop n\sim j 2^k}\sum\limits_{m\in I\cap\mathbb{Z},\atop m\sim (j+1) 2^k}u_m u_n\bigg\|_{L^2([0,T]\times \mathbb{R})}^2\bigg)^{1/2}.
\end{align}

Denote $u_{j,k}:=\sum\limits_{n\in I,\atop n\sim j2^k}u_n$, from proposition \ref{UVprop2} we can write as a sum $u_{j,k}=u_{1,j,k}+u_{2,j,k}$ with the estimate
\begin{align}\label{L40}
\frac{1}{|I|^{\beta}}\|u_{1,j,k}\|_{U_{A}^2}+e^{\epsilon |I|^{\beta}}\| u_{2,j,k}\|_{U_{A}^8}\lesssim \| u_{j,k}\|_{V_{A}^2}.
\end{align}

Then the estimate \eqref{lebesgue4b} will be continued by four terms. For the term containing $u_{1,j,k}$ and $u_{1,j+1,k}$, which will be denoted as $I_1$.
\begin{align}\label{L41}
I_1&\lesssim \sum_{k\in \mathbb{N}}\bigg(\sum_{j\in \mathbb{N}^+} \|u_{1,j,k}u_{1,j+1,k}\|^{2\theta}_{L^2}\|u_{1,j,k}u_{1,j+1,k}\|^{2(1-\theta)}_{L^2([0,T]\times \mathbb{R})}\bigg)^{1/2}\nonumber\\
&\lesssim \sum_{k\in \mathbb{N}}\bigg(\sum_{j\in \mathbb{N}^+} \|u_{1,j,k}u_{1,j+1,k}\|^{2\theta}_{L^2}\|u_{1,j,k}\|^{2(1-\theta)}_{L^4([0,T]\times \mathbb{R})}\|u_{1,j+1,k}\|^{2(1-\theta)}_{L^4([0,T]\times \mathbb{R})}\bigg)^{1/2}\nonumber\\
&\lesssim T^{(1-\theta)/4}\sum_{k\in \mathbb{N}}\bigg(\sum_{j\in \mathbb{N}^+} \|u_{1,j,k}u_{1,j+1,k}\|^{2\theta}_{L^2}\|u_{1,j,k}\|^{2(1-\theta)}_{L^8_{t\in[0,T]}L^4_x}\|u_{1,j+1,k}\|^{2(1-\theta)}
_{L^8_{t\in[0,T]}L^4_x}\bigg)^{1/2}.
\end{align}
Since $|m-n|\sim 2^k$, $|m+n|\sim (2j+1)2^k\gtrsim\sqrt{j(j+1)}2^k$, we have the bilinear estimates
\begin{align}
 \|u_{1,j,k}u_{1,j+1,k}\|_{L^2} \lesssim (j(j+1))^{-1/4}2^{-k}\|u_{1,j,k}\|_{U^2_A} \|u_{1,j+1,k}\|_{U^2_A}.
\end{align}
Combining with Strichartz estimate, \eqref{L41} is dominated by
\begin{align*}
&\lesssim  T^{(1-\theta)/4}\sum_{k\in \mathbb{N}}\bigg(\sum_{j\in \mathbb{N}^+} (j(j+1))^{-\theta/2}2^{-2k\theta} \|u_{1,j,k}\|^{2\theta}_{U^2_A}\|u_{1,j+1,k}\|^{2\theta}_{U^2_A}\nonumber\\
&\quad\quad\quad(j2^k)^{-(1-\theta)/4}\|u_{1,j,k}\|^{2(1-\theta)}_{U^8_A}((j+1)2^k)^{-(1-\theta)/4}\|u_{1,j+1,k}\|
^{2(1-\theta)}_{U^8_A}\bigg)^{1/2}\nonumber\\
&\lesssim T^{(1-\theta)/4}\sum_{k\in \mathbb{N}} \bigg(\sum_{j\in \mathbb{N}^+} (j2^k)^{-(1+\theta)/4}((j+1)2^k)^{-(1+\theta)/4}2^{-k\theta}\|u_{1,j,k}\|^2_{U^2_A}\|u_{1,j+1,k}\|^2_{U^2_A}
\bigg)^{1/2}.
\end{align*}
By applying \eqref{L40} and the orthogonality in $V^2$, it follows that
\begin{align}\label{L43}
&\lesssim T^{(1-\theta)/4}|I|^{2\beta}\sum_{k\in \mathbb{N}}\bigg(\sum_{j\in \mathbb{N}^+} (j2^k)^{-1/4}((j+1)2^k)^{-1/4}2^{-k\theta}\|u_{1,j,k}\|^2_{V^2_A}\|u_{1,j+1,k}\|^2_{V^2_A}
\bigg)^{1/2}\nonumber\\
&\lesssim T^{(1-\theta)/4}|I|^{2\beta}\sum_{k\in \mathbb{N}}\bigg(\sum_{j\in \mathbb{N}^+} (j2^k)^{-1/4}((j+1)2^k)^{-1/4}2^{-k\theta}\nonumber\\
&\quad\quad\quad\quad\quad\quad\cdot\bigg(\sum\limits_{n\in I,\atop n\sim j2^k}\|u_n\|^2_{V^2_A}\bigg)\bigg(\sum\limits_{m\in I,\atop m\sim (j+1)2^k}\|u_m\|^2_{V^2_A}\bigg)
\bigg)^{1/2}\nonumber\\
&\lesssim T^{(1-\theta)/4}|I|^{2\beta}\sum_{k\in \mathbb{N}}\Bigg(\sum_{j\in \mathbb{N}^+} (j2^k)^{-1/4}((j+1)2^k)^{-1/4}2^{-k\theta}2^k\nonumber\\
&\quad\quad\quad\quad\quad\quad\cdot\bigg(\sum\limits_{n\in I,\atop n\sim j2^k}\|u_n\|^4_{V^2_A}\bigg)^{1/2}\bigg(\sum\limits_{m\in I,\atop m\sim (j+1)2^k}\|u_m\|^4_{V^2_A}\bigg)^{1/2}
\Bigg)^{1/2}\nonumber\\
&\lesssim T^{(1-\theta)/4}|I|^{2\beta+(1-\theta)/2}\ln|I|\|u\|^2_{X^{-1/8}_{4,A}}\nonumber\\
&\lesssim T^{(1-\theta)/4}|I|^{2\beta+(1-\theta)/2}\|u\|^2_{X^{-1/8}_{4,A}},
\end{align}
where the last inequality is by using $2^{k(1-\theta)/2}\lesssim |I|^{(1-\theta)/2}$ and H$\ddot{{\rm o}}$lder's inequality. For the rest three terms we will do in a uniform way. We take the term containing $u_{2,j,k}$ and $u_{2,j+1,k}$ for example, and denote it as $I_2$.
\begin{align}\label{L44}
I_2&\lesssim \sum_{k\in \mathbb{N}}\bigg(\sum_{j\in \mathbb{N}^+} \|u_{2,j,k}u_{2,j+1,k}\|^{2}_{L^2([0,T]\times \mathbb{R})}\bigg)^{1/2}\nonumber\\
&\lesssim T^{1/4}\sum_{k\in \mathbb{N}}\bigg(\sum_{j\in \mathbb{N}^+} \|u_{2,j,k}\|^{2}_{L^8_{t\in[0,T]}L^4_x}\|u_{2,j+1,k}\|^{2}_{L^8_{t\in[0,T]}L^4_x}\bigg)^{1/2}\nonumber\\
&\lesssim T^{1/4}\sum_{k\in \mathbb{N}}\bigg(\sum_{j\in \mathbb{N}^+} (j2^k)^{-1/4}\|u_{2,j,k}\|^{2}_{U^8_A}((j+1)2^k)^{-1/4}\|u_{2,j+1,k}\|^{2}_{U^8_A}\bigg)^{1/2}.
\end{align}
By applying \eqref{L40} and the orthogonality in $V^2$ again, \eqref{L44} follows that
\begin{align}
&\lesssim T^{1/4}e^{-2\epsilon|I|^\beta}\sum_{k\in \mathbb{N}}\bigg(\sum_{j\in \mathbb{N}^+} (j2^k)^{-1/4}\|u_{j,k}\|^2_{V^2_A}((j+1)2^k)^{-1/4}2^{-k\theta}\|u_{j+1,k}\|^2_{V^2_A}
\bigg)^{1/2}\nonumber\\
&\lesssim T^{1/4}e^{-2\epsilon|I|^\beta}\sum_{k\in \mathbb{N}}\bigg(\sum_{j\in \mathbb{N}^+} (j2^k)^{-1/4}\bigg(\sum\limits_{n\in I,\atop n\sim j2^k}\|u_n\|^2_{V^2_A}\bigg)((j+1)2^k)^{-1/4}\bigg(\sum\limits_{m\in I,\atop m\sim (j+1)2^k}\|u_m\|^2_{V^2_A}\bigg)
\bigg)^{1/2}\nonumber\\
&\lesssim T^{1/4}e^{-2\epsilon|I|^\beta}\sum_{k\in \mathbb{N}}\Bigg(\sum_{j\in \mathbb{N}^+} (j2^k)^{-1/4}((j+1)2^k)^{-1/4}2^k\bigg(\sum\limits_{n\in I,\atop n\sim j2^k}\|u_n\|^4_{V^2_A}\bigg)^{1/2}\nonumber\\
&\quad\quad\quad\quad\quad\quad\cdot\bigg(\sum\limits_{m\in I,\atop m\sim (j+1)2^k}\|u_m\|^4_{V^2_A}\bigg)^{1/2}
\Bigg)^{1/2}\nonumber\\
&\lesssim T^{1/4}e^{-2\epsilon|I|^\beta}|I|^{1/2}\ln|I|\|u\|^2_{X^{-1/8}_{4,A}}\nonumber\\
&\lesssim T^{1/4}\|u\|^2_{X^{-1/8}_{4,A}}.
\end{align}
Thus we complete the proof of \eqref{lebesgue4}. In particular, for $1 \lesssim  |I| <\infty$ and $0<T<1$, taking $\beta$ and $1-\theta$ small sufficiently, we have
\begin{align}
\|u_I\|_{L^4_{x,t\in [0,T]}} \lesssim   T^{\varepsilon/4}|I|^{\varepsilon}  \|u\|_{X^{-1/8}_{4,A}(I)}.
\end{align}
In the end we can obtain \eqref{lebesgue4a} by H$\ddot{{\rm o}}$lder inequality.$\hfill\Box$
\begin{lem}\label{V2toX}
Let $I\subset \mathbb{R}$ with $1\lesssim  |I|<\infty$, $2\leq q \leq\infty$, we have
\begin{align}
\|u_I\|_{L^\infty_t L^2_{x} \cap V^2_A} & \lesssim    |I|^{1/2-1/q} \max_{\lambda \in I} \langle \lambda\rangle^{-1/4} \|u\|_{X^{1/4}_{q,A}(I)}.  \label{V2}
\end{align}
\end{lem}
{\bf Proof.} Using $V^2_A \subset L^\infty_t L^2_x$, the orthogonality in $V^2$ and H\"older's inequality one by one, we have
\begin{align}
\|u_I\|_{L^\infty_t L^2_{x} \cap V^2_A} &\lesssim  \bigg(\sum_{\lambda\in I}\|u_\lambda\|^2_{V^2_A}\bigg)^{1/2}
  \lesssim\max_{\lambda\in I}\langle\lambda\rangle^{-1/4}\bigg(\sum_{\lambda\in I}\big(\langle\lambda\rangle^{1/4}\|u_\lambda\|\big)^2_{V^2_A}\bigg)^{1/2}\nonumber\\
  &\lesssim|I|^{1/2-1/q} \max_{\lambda \in I} \langle \lambda\rangle^{-1/4} \|u\|_{X^{1/4}_{q,A}(I)}.
\end{align}

\section{Trilinear estimates}

At first, we apply the duality  to the norm calculation (Proposition \ref{UVprop4}) to the inhomogeneous part of the solution of mKdV in $X^s_{q,A}$. It is known that \eqref{mKdV} is equivalent to the following integral equation:
\begin{align}
u(x,t)=e^{-t\partial_x^3} u_0-\mathcal{A}((u^3)_x) ,   \label{ImKdV}
\end{align}
where
$$
e^{-t\partial_x^3}=\mathscr{F}^{-1}e^{{\rm i}t\xi^3}\mathscr{F},\ \ \ \mathcal{A} (f) = \int^t_0 e^{-(t-\tau)\partial_x^3} f(\tau) d\tau.
$$
By Propositions \ref{UVprop3} and \ref{UVprop4}, we see that, for ${\rm supp} \ v \subset \mathbb{R}\times [0,T]$, $1\leq q<\infty$,
\begin{align}
\|\mathcal{A} (f) \|_{X^{1/4}_{q,A}} &  =\|e^{t\partial_x^3}\mathcal{A} (f) \|_{X^{1/4}_{q}}\nonumber\\
&=\sup \left\{ \left|B\left(\int^t_0  e^{\tau\partial_x^3} f(\tau)d\tau, v \right) \right| : \  \|v\|_{Y^{-1/4}_{q'}} \leq 1 \right\} \nonumber\\
& \leq \sup_{\|v\|_{Y^{-1/4}_{q'}} \leq 1} \left|\int_{ [0,T]} \langle  e^{t\partial_x^3}f (t), \  v(t) \rangle dt  \right| \nonumber\\
& \leq \sup_{\|v\|_{Y^{-1/4}_{q'}} \leq 1} \left|\int_{ [0,T]} \langle f (t), \  e^{-t\partial_x^3} v(t) \rangle dt  \right| \nonumber\\
& \leq \sup_{\|v\|_{Y^{-1/4}_{q', A}} \leq 1} \left|\int_{ [0,T]} \langle f (t), \  v(t) \rangle dt  \right|.  \label{normest1}
\end{align}
For $q=\infty$, we have
\begin{align}
\|\mathcal{A} (f) \|_{X^{1/4}_{\infty,A}} &  =\|e^{t\partial_x^3}\mathcal{A} (f) \|_{X^{1/4}_{\infty}}\nonumber\\
&=\sup_{\lambda\in \mathbb{Z}}\ \langle\lambda\rangle^{1/4}\bigg\|\Box_\lambda\int^t_0  e^{\tau\partial_x^3} f(\tau)d\tau\bigg\|_{U^2}\nonumber\\
&\leq\sup_{\lambda\in \mathbb{Z}}\ \langle\lambda\rangle^{1/4}\sup_{\|v^{(\lambda)}\|_{V^2} \leq 1} \left|\int_{ [0,T]} \langle \Box_{\lambda} e^{t\partial_x^3}f (t), \  v^{(\lambda)}(t) \rangle dt  \right| \nonumber\\
&\leq\sup_{\lambda\in \mathbb{Z}}\ \langle\lambda\rangle^{1/4}\sup_{\|v^{(\lambda)}\|_{V^2} \leq 1} \left|\int_{ [0,T]} \langle f (t), \  \Box_{\lambda} e^{-t\partial_x^3}v^{(\lambda)}(t) \rangle dt  \right| \nonumber\\
&\leq\sup_{\lambda\in \mathbb{Z}}\ \langle\lambda\rangle^{1/4}\sup_{\|v^{(\lambda)}\|_{V_A^2} \leq 1} \left|\int_{ [0,T]} \langle f (t), \  \Box_{\lambda} v^{(\lambda)}(t) \rangle dt  \right| .  \label{normest2}
\end{align}

To prove Theorem \ref{theorem1}, we need to control the second term of the integral equation \eqref{ImKdV} in $X^{1/4}_{q,A} (2\leq q\leq \infty)$. More precisely, we want to prove the following lemma.
\begin{lem}
For $2\leq q\leq \infty$, there exists $\varepsilon>0$ such that
\begin{align}
 \left\|\int^t_0 e^{-(t-\tau)\partial_x^3} (u^3)_x(\tau)\, d\tau  \right\|_{X^{1/4}_{q,A}}\lesssim T^{\varepsilon}\|u\|^{3}_{X^{1/4}_{q,A}}.  \label{trilinear1}
\end{align}
\end{lem}
{\bf Proof.} When $2\leq q< \infty$, in view of \eqref{normest1},  it suffices to show that
\begin{align}
&  \left|\int_{\mathbb{R}\times [0,T]} \overline{v} u^2\partial_x u \, dxdt \right|\lesssim T^{\varepsilon}\|u\|^{3}_{X^{1/4}_{q,A}}  \|v\|_{Y^{-1/4}_{q',A}} . \label{trilinear2}
\end{align}
We perform a uniform decomposition with $u,v$ in the left hand side of \eqref{trilinear2}, it suffices to prove that
\begin{align}
\left\lvert \sum_{\lambda_0,...,\lambda_3} \langle \lambda_0 \rangle^{1/4}\int_{[0,T]\times\mathbb{R}} \overline{v}_{\lambda_0} u_{\lambda_1}u_{\lambda_2}\partial_{x}u_{\lambda_3} \, dxdt \right \rvert \lesssim T^{\varepsilon} \| u\| ^{3}_{X^{1/4}_{q,A}}  \|v\|_{Y^{0}_{q',A}}. \label{trilinear3}
\end{align}
When $q=\infty$, in view of \eqref{normest2}, it suffices to show that, for any fixed $\lambda\in\mathbb{Z}$,
\begin{align}
\left\lvert \sum_{\lambda_1,\lambda_2,\lambda_3} \langle \lambda\rangle^{1/4}\int_{[0,T]\times\mathbb{R}} \overline{\Box_\lambda v^{(\lambda)}} u_{\lambda_1}u_{\lambda_2}\partial_{x}u_{\lambda_3} \, dxdt \right \rvert \lesssim T^{\varepsilon} \| u\| ^{3}_{X^{1/4}_{\infty,A}}  \|v^{(\lambda)}\|_{V^2_A}. \label{trilinear4}
\end{align}
\subsection {$q=\infty$, Proof of \eqref{trilinear4}.}

For convenience, denote $\lambda$ as $\lambda_0$, $\Box_\lambda v^{(\lambda)}=:v_\lambda=v_{\lambda_0}$.
In order to keep the left hand side of \eqref{trilinear4} nonzero, we have the frequency constraint condition (FCC)
\begin{equation}\label{FCC}
\lambda_1+\lambda_2+\lambda_3\approx\lambda_0
\end{equation}
and dispersion modulation constraint condition (DMCC)
\begin{equation}\label{DMCC}
\max_{0\leq k \leq 3} |\xi^3_k - \tau_k|  \gtrsim \bigg|(\xi^3_0 - \tau_0)-\sum_{1\leq k\leq 3}(\xi^3_k - \tau_k)\bigg|\gtrsim|(\xi_0- \xi_1)(\xi_0-\xi_2)(\xi_0-\xi_3)|.
\end{equation}
It suffices to consider the cases that $\lambda_0$ is maximal or secondly maximal number in $\lambda_0,...,\lambda_3$ (In the opposite case, one can replace $\lambda_0,...,\lambda_3$ with $-\lambda_0,...,-\lambda_3$).

{\bf Step 1.} We assume that $\lambda_0=\max_{0\leq k\leq 3}\lambda_k$. From the frequency constraint condition (FCC) $\lambda_0\approx\lambda_1+\lambda_2+\lambda_3$, we know that the non-trivial case is that $\lambda_0\gg 0$ ( The case $\lambda_0\ll 0$ never happens due to the condition (FCC). In addition, the case $|\lambda_0|\lesssim 1$, which leads to $\max_{0\leq k\leq 3}|\lambda_k|\lesssim 1$, implies that the summation in \eqref{trilinear4} has at most finite terms). Furthermore, in view of $\lambda_0=\max_{0\leq k\leq 3}\lambda_k$, $\lambda_0\approx\lambda_1+\lambda_2+\lambda_3$, and $\lambda_0\gg 0$, we see that $\lambda_0\approx\max_{0\leq k\leq 3}|\lambda_k|\gg0$. For convenience, we can take
\begin{align}\label{step1a}
\lambda_0=\max_{0\leq k\leq 3}|\lambda_k|\gg0.
\end{align}
By the symmetry, we can assume $\lambda_1\geq\lambda_2$. Then $\lambda_0,...,\lambda_3$ have the following three orders:
\begin{align*}
{\rm Order\,1:}\quad\quad\lambda_0\geq \lambda_3 \geq \lambda_1 \geq \lambda_2\ ;\\
{\rm Order\,2:}\quad\quad\lambda_0\geq \lambda_1 \geq \lambda_3 \geq \lambda_2\ ;\\
{\rm Order\,3:}\quad\quad\lambda_0\geq \lambda_1 \geq \lambda_2 \geq \lambda_3\ .
\end{align*}
We just take Order 1 for example because the other two orders are similar and even more easier (noticing that the derivative located in $u_{\lambda_3}$).

{\bf Order 1}: $\lambda_0\geq \lambda_3 \geq \lambda_1 \geq \lambda_2$.  For short, considering the higher and lower frequency of $\lambda_k$,  we use the following notations:
$$
\left\{
\begin{array}{l}
\lambda_k \in h \Leftrightarrow \lambda_k\in [3\lambda_0/4, \lambda_0]; \\
\lambda_k \in h_- \Leftrightarrow \lambda_k\in [-\lambda_0, -3\lambda_0/4];\\
\lambda_k \in l \Leftrightarrow \lambda_k\in [0, 3 \lambda_0/4];\\
\lambda_k \in l_- \Leftrightarrow \lambda_k\in [-3\lambda_0/4, 0].
\end{array}
\right.
$$
Then we divide Order 1 into several cases.

{\bf Case 1}: $\lambda_3\in h$ and $\lambda_1\in h$.  In consideration of (FCC),  we easily see that $\lambda_0, \lambda_1, \lambda_2,\lambda_3$ satisfy the following frequency constraint condition:
\begin{align}\label{FCCl}
\lambda_0=\lambda_1+\lambda_2+\lambda_3+l, \ \ \ |l|\leq 10.
\end{align}
We know that this case implies that $\lambda_2\in [-\lambda_0,-\lambda_0/2-l]$. We do dyadic decomposition for $u_{\lambda_1}$, $u_{\lambda_2}$ and $u_{\lambda_3}$, and keep using uniform decomposition for $v_{\lambda_0}$. Let us denote $I_0=[0,1)$, $I_j=[2^{j-1}, 2^j)$, $j\geq 1$.  We decompose $\lambda_1, \lambda_2, \lambda_3$ by:
\begin{align}
\lambda_k \in [3\lambda_0/4, \lambda_0]= \bigcup_{j_k\geq 0} \lambda_0-I_{j_k}, \ k=1,3;\ \   \lambda_2 \in [-\lambda_0, -\lambda_0/2-l]= \bigcup_{j_2\geq 0}-\lambda_0+I_{j_2}.    \nonumber
\end{align}
From $\lambda_3\geq \lambda_1$ we know that $j_3\leq j_1$. In view of condition (FCC), we see that $j_2\approx j_1$. It follows that
\begin{align*}
0\leq j_3\leq j_1\approx j_2\leq \log_2{\lambda_0}.
\end{align*}
In the following discussion, we shall omit the condition $j_k\in [0, \log_2{\lambda_0}]$, $k=1,2,3$, for convenience, but it is always satisfied in Step 1.
We denote the left hand side of \eqref{trilinear4} as $\mathscr{L}_{hhhh_-}(u, v)$, and divide it into three parts:
\begin{align*}
 \mathscr{L}_{hhhh_-}(u, v):&=\sum_{j_3\leq j_1\approx j_2} \langle \lambda_0\rangle^{1/4}\int_{[0,T]\times\mathbb{R}} |\overline{v_{\lambda_0}} u_{\lambda_0-I_{j_1}}u_{-\lambda_0+I_{j_2}}\partial_{x}u_{\lambda_0-I_{j_3}} |\, dxdt\nonumber\\
 &=\bigg(\sum_{j_3\leq j_1\approx j_2\lesssim1}+\sum_{j_3\lesssim 1\ll j_1\approx j_2}+\sum_{1\ll j_3\leq j_1\approx j_2}\bigg)\nonumber\\
 &\quad\quad\langle \lambda_0\rangle^{1/4}\int_{[0,T]\times\mathbb{R}} |\overline{v_{\lambda_0}} u_{\lambda_0-I_{j_1}}u_{-\lambda_0+I_{j_2}}\partial_{x}u_{\lambda_0-I_{j_3}} |\, dxdt\nonumber\\
 &=:\mathscr{L}^l_{hhhh_-}(u, v)+\mathscr{L}^m_{hhhh_-}(u, v)+\mathscr{L}^h_{hhhh_-}(u, v).
\end{align*}
It is easy to see that in $\mathscr{L}^l_{hhhh_-}(u, v)$, $\lambda_0\approx\lambda_3\approx\lambda_1\approx -\lambda_2$ holds. Therefore, by H$\ddot{\rm o}$lder inequality and Strichartz estimate, we have
\begin{align}
 \mathscr{L}^l_{hhhh_{-}}(u, v)&\lesssim\langle \lambda_0\rangle^{5/4} \|\overline{v_{\lambda_0}}\|_{L^4_{x,t}}\|u_{\lambda_0}\|^2_{L^4_{x,t}}\|u_{-\lambda_0}\|_{L^4_{x,t}}\nonumber\\
&\lesssim  T^{1/2}\langle \lambda_0\rangle^{5/4} \|\overline{v_{\lambda_0}}\|_{L^8_tL^4_x}\|u_{\lambda_0}\|^2_{L^8_tL^4_x}\|u_{-\lambda_0}\|_{L^8_tL^4_x}\nonumber\\
&\lesssim T^{1/2}\langle \lambda_0\rangle^{3/4} \|\overline{v_{\lambda_0}}\|_{U^8_A}\|u_{\lambda_0}\|^2_{U^8_A}\|u_{-\lambda_0}\|_{U^8_A}\nonumber\\
&\lesssim  T^{1/2}\langle \lambda_0\rangle^{3/4} \|\overline{v_{\lambda_0}}\|_{V^2_A}\|u_{\lambda_0}\|^2_{U^2_A}\|u_{-\lambda_0}\|_{U^2_A}\nonumber\\
&\lesssim T^{1/2} \|v^{(\lambda_0)}\|_{V^2_A}\|u\|^{3}_{X^{1/4}_{\infty,A}}. \nonumber
\end{align}
In $\mathscr{L}^m_{hhhh_-}(u, v)$, we see that the frequency of $v_{\lambda_0}$ and $u_{\lambda_0-I_{j_3}}$ are localized near $\lambda_0$, which are far away from the frequency of $u_{\lambda_0-I_{j_1}}$ and the reciprocal frequency of  $u_{-\lambda_0+I_{j_2}}$. Thus we can use bilinear estimate \eqref{bilinear3}, Lemma \ref{V2toX}, and H\"older's inequality to obtain that
\begin{align*}
\mathscr{L}^m_{hhhh_-}(u, v)&\lesssim\sum_{j_3\lesssim 1\ll j_1\approx j_2}\langle \lambda_0\rangle^{1/4} \|\overline{v_{\lambda_0}}u_{\lambda_0-I_{j_1}}\|_{L^2_{x,t}}\|u_{-\lambda_0+I_{j_2}}
\partial_{x}u_{\lambda_0-I_{j_3}}\|_{L^2_{x,t}}\nonumber\\
&\lesssim  T^{\varepsilon/2}\sum_{j_3\lesssim 1\ll j_1\approx j_2}\langle \lambda_0\rangle^{5/4}\langle \lambda_0\rangle^{-1/2+\varepsilon}(2^{j_1})^{-1/2+\varepsilon}\|v_{\lambda_0}\|_{V^2_A} \|u_{\lambda_0-I_{j_1}}\|_{V^2_A}\nonumber\\
&\quad\quad\quad\quad\times \langle \lambda_0\rangle^{-1/2+\varepsilon}(2^{j_2})^{-1/2+\varepsilon}\| u_{-\lambda_0+I_{j_2}}\|_{V^2_A}\|u_{\lambda_0-I_{j_3}}\|_{V^2_A}\nonumber\\
&\lesssim  T^{\varepsilon/2}\sum_{j_3\lesssim 1\ll j_1\approx j_2}\langle \lambda_0\rangle^{1/4+2\varepsilon}(2^{j_1})^{-1/2+\varepsilon}(2^{j_2})^{-1/2+\varepsilon}\|v_{\lambda_0}\|_{V^2_A}\nonumber\\
&\quad\quad\quad\quad\times(2^{j_1})^{1/2}(2^{j_2})^{1/2}(2^{j_3})^{1/2}\langle\lambda_0\rangle^{-3/4}\| u\|^3_{X^{1/4}_{\infty,A}}\nonumber\\
&\lesssim  T^{\varepsilon/2}\langle \lambda_0\rangle^{-1/2+4\varepsilon}
\|v_{\lambda_0}\|_{V^2_A}\|u\|^3_{X^{1/4}_{\infty,A}}\nonumber\\
&\lesssim  T^{\varepsilon/2}\|v_{\lambda_0}\|_{V^2_A}\|u\|^3_{X^{1/4}_{\infty,A}},
\end{align*}
where the last inequality is obtained by taking $\varepsilon\leq 1/8$.

Now we estimate $\mathscr{L}^h_{hhhh_-}(u, v)$. In view of (DMCC) \eqref{DMCC}, we have the highest dispersion modulation satisfying
\begin{align*}
\max_{0\leq k \leq 3} |\xi^3_k - \tau_k|  \gtrsim \langle\lambda_0\rangle\cdot 2^{j_1}\cdot 2^{j_3}.
\end{align*}

If $v_{\lambda_0}$ has the highest dispersion modulation, we divide $\mathscr{L}^h_{hhhh_-}(u, v)$ into two parts.
\begin{align}\label{hhhh-}
&\mathscr{L}^{h1}_{hhhh_-}(u, v)+\mathscr{L}^{h2}_{hhhh_-}(u, v)\nonumber\\
:=&\bigg(\sum_{1\ll j_3\ll j_1\approx j_2}+\sum_{1\ll j_3\approx j_1\approx j_2}\bigg)\langle \lambda_0\rangle^{1/4}\int_{[0,T]\times\mathbb{R}} |\overline{v_{\lambda_0}} u_{\lambda_0-I_{j_1}}u_{-\lambda_0+I_{j_2}}\partial_{x}u_{\lambda_0-I_{j_3}} |\, dxdt.
\end{align}
For $\mathscr{L}^{h1}_{hhhh_-}(u, v)$, we can use the bilinear estimate due to $j_2\gg j_3$. By H\"older's inequality we have
\begin{align*}
\mathscr{L}^{h1}_{hhhh_-}(u, v)\lesssim \sum_{1\ll j_3\ll j_1\approx j_2}\langle \lambda_0\rangle^{1/4}\|\overline{v_{\lambda_0}}\|_{L^2_tL^{\infty}_x} \|u_{\lambda_0-I_{j_1}}\|_{L^\infty_tL^{2}_x}\| u_{-\lambda_0+I_{j_2}}\partial_{x}u_{\lambda_0-I_{j_3}}\|_{L^2_{x,t}}.
\end{align*}
Using $\|v_{\lambda_0}\|_{L^\infty_x}\lesssim \|v_{\lambda_0}\|_{L^2_x}$, $V^2_A\subset L^\infty_t L^2_x$, the dispersion modulation decay \eqref{dispersiondecay}, the bilinear estimate \eqref{bilinear3} and Lemma \ref{V2toX}, we have
\begin{align}\label{hhhh-1}
\mathscr{L}^{h1}_{hhhh_-}(u, v)&\lesssim \sum_{1\ll j_3\ll j_1\approx j_2} \langle \lambda_0\rangle^{5/4}\langle \lambda_0\rangle^{-1/2}(2^{j_1})^{-1/2}(2^{j_3})^{-1/2}\|v_{\lambda_0}\|_{V^2_A} \|u_{\lambda_0-I_{j_1}}\|_{V^2_A}\nonumber\\
&\quad\quad \times T^{\varepsilon/4}\langle \lambda_0\rangle^{-1/2+\varepsilon}(2^{j_2})^{-1/2+\varepsilon}\| u_{-\lambda_0+I_{j_2}}\|_{V^2_A}\|u_{\lambda_0-I_{j_3}}\|_{V^2_A}\nonumber\\
&\lesssim  T^{\varepsilon/4}\sum_{1\ll j_3\ll j_1\approx j_2}\langle \lambda_0\rangle^{1/4+\varepsilon}(2^{j_1})^{-1/2}(2^{j_2})^{-1/2+\varepsilon}(2^{j_3})^{-1/2}\|v_{\lambda_0}\|_{V^2_A}\nonumber\\
&\quad\quad\quad\quad\times(2^{j_1})^{1/2}(2^{j_2})^{1/2}(2^{j_3})^{1/2}\langle\lambda_0\rangle^{-3/4}\| u\|^3_{X^{1/4}_{\infty,A}}\nonumber\\
&\lesssim  T^{\varepsilon/4}\langle \lambda_0\rangle^{-1/2+\varepsilon}\bigg(\sum_{j_2}(2^{j_2})^{\varepsilon}\cdot j_2\bigg)\|v_{\lambda_0}\|_{V^2_A}\| u\|^3_{X^{1/4}_{\infty,A}}.
\end{align}
Noticing that $2^{j_2}\lesssim \langle \lambda_0\rangle$ and $j_2\lesssim (2^{j_2})^\varepsilon$, we can take $\varepsilon\leq 1/6$ such that
\begin{align*}
\mathscr{L}^{h1}_{hhhh_-}(u, v)&\lesssim  T^{\varepsilon/4}\langle \lambda_0\rangle^{-1/2+3\varepsilon}\|v_{\lambda_0}\|_{V^2_A}\| u\|^3_{X^{1/4}_{\infty,A}}\nonumber\\
&\lesssim  T^{\varepsilon/4}\|v_{\lambda_0}\|_{V^2_A}\| u\|^3_{X^{1/4}_{\infty,A}}.
\end{align*}
For $\mathscr{L}^{h2}_{hhhh_-}(u, v)$, by H\"older's inequality, $\|v_{\lambda_0}\|_{L^\infty_x}\lesssim \|v_{\lambda_0}\|_{L^2_x}$, $V^2_A\subset L^\infty_t L^2_x$, the dispersion modulation decay \eqref{dispersiondecay}, the $L^4$ estimate \eqref{lebesgue4a}  and Lemma \ref{V2toX}, we have
\begin{align}\label{hhhh-2}
\mathscr{L}^{h2}_{hhhh_-}(u, v)&\lesssim \sum_{1\ll j_3\approx j_1\approx j_2}\langle \lambda_0\rangle^{5/4}\|\overline{v_{\lambda_0}}\|_{L^2_tL^{\infty}_x} \|u_{\lambda_0-I_{j_1}}\|_{L^\infty_tL^{2}_x}\| u_{-\lambda_0+I_{j_2}}\|_{L^4_{x,t}}\|u_{\lambda_0-I_{j_3}}\|_{L^4_{x,t}}\nonumber\\
&\lesssim \sum_{1\ll j_3\approx j_1\approx j_2}\langle \lambda_0\rangle^{5/4}\langle \lambda_0\rangle^{-1/2}(2^{j_1})^{-1/2}(2^{j_3})^{-1/2}\|v_{\lambda_0}\|_{V^2_A} \|u_{\lambda_0-I_{j_1}}\|_{V^2_A}\nonumber\\
&\quad\quad \times T^{\varepsilon/2}\langle \lambda_0\rangle^{-3/4}(2^{j_2})^{1/4+\varepsilon}(2^{j_3})^{1/4+\varepsilon}\| u\|^2_{X^{1/4}_{\infty,A}}\nonumber\\
&\lesssim T^{\varepsilon/2}\sum_{j_3\approx j_2}\langle \lambda_0\rangle^{-1/4}(2^{j_2})^{1/4+\varepsilon}(2^{j_3})^{-1/4+\varepsilon}
\|v_{\lambda_0}\|_{V^2_A}\|u\|^3_{X^{1/4}_{\infty,A}}\nonumber\\
&\lesssim T^{\varepsilon/2}\langle \lambda_0\rangle^{-1/4+2\varepsilon}\|v_{\lambda_0}\|_{V^2_A}\|u\|^3_{X^{1/4}_{\infty,A}}\nonumber\\
&\lesssim T^{\varepsilon/2}\|v_{\lambda_0}\|_{V^2_A}\|u\|^3_{X^{1/4}_{\infty,A}},
\end{align}
where the last inequality is by taking $\varepsilon\leq 1/8$.

If $u_{\lambda_0-I_{j_1}}$ has the highest dispersion modulation, we just take $L^\infty_{x,t}$ and $L^2_{x,t}$ norms to $v_{\lambda_0}$ and $u_{\lambda_0-I_{j_1}}$, respectively, then
\begin{align}
\|\overline{v_{\lambda_0}}\|_{L^\infty_{x,t}} \|u_{\lambda_0-I_{j_1}}\|_{L^2_{x,t}}\lesssim \|v_{\lambda_0}\|_{V^2_A} \langle \lambda_0\rangle^{-1/2}(2^{j_1})^{-1/2}(2^{j_3})^{-1/2}\|u_{\lambda_0-I_{j_1}}\|_{V^2_A},
\end{align}
where we use the fact $\|v_{\lambda_0}\|_{L^\infty_{x,t}}\lesssim \|v_{\lambda_0}\|_{L^\infty_tL^2_x}\lesssim \|v_{\lambda_0}\|_{V^2_A}$ and the dispersion modulation decay \eqref{dispersiondecay} to $u_{\lambda_0-I_{j_1}}$. Then this case reduces to the same estimate as that when $v_{\lambda_0}$ has the highest dispersion modulation.

If $u_{-\lambda_0+I_{j_2}}$ has the highest dispersion modulation, we still divide $\mathscr{L}^h_{hhhh_-}(u, v)$ into two parts as \eqref{hhhh-}.
For $\mathscr{L}^{h1}_{hhhh_-}(u, v)$, we can take $L^\infty_{x,t}$, $L^2_{x,t}$ and $L^2_{x,t}$ norms to $v_{\lambda_0}$, $u_{-\lambda_0+I_{j_2}}$ and $u_{\lambda_0-I_{j_1}}\partial_xu_{\lambda_0-I_{j_3}}$, respectively. By H\"older's inequality, the dispersion modulation decay \eqref{dispersiondecay}, the bilinear estimate \eqref{bilinear3} and Lemma \ref{V2toX}, we have
\begin{align*}
\mathscr{L}^{h1}_{hhhh_-}(u, v)&\lesssim \sum_{1\ll j_3\ll j_1\approx j_2}\langle \lambda_0\rangle^{1/4}\|\overline{v_{\lambda_0}}\|_{L^\infty_{x,t}} \|u_{-\lambda_0+I_{j_2}}\|_{L^2_{x,t}}\| u_{\lambda_0-I_{j_1}}\partial_{x}u_{\lambda_0-I_{j_3}}\|_{L^2_{x,t}}\nonumber\\
&\lesssim \sum_{1\ll j_3\ll j_1\approx j_2} \langle \lambda_0\rangle^{5/4}\langle \lambda_0\rangle^{-1/2}(2^{j_1})^{-1/2}(2^{j_3})^{-1/2}\|v_{\lambda_0}\|_{V^2_A}\| u_{-\lambda_0+I_{j_2}}\|_{V^2_A} \nonumber\\
&\quad\quad \times T^{\varepsilon/4}\langle \lambda_0\rangle^{-1/2+\varepsilon}(2^{j_1})^{-1/2+\varepsilon}\|u_{\lambda_0-I_{j_1}}\|_{V^2_A}\|u_{\lambda_0-I_{j_3}}\|_{V^2_A},
\end{align*}
which is the same as the right hand side of the first inequality in \eqref{hhhh-1} (noticing that $j_1\approx j_2$).

For $\mathscr{L}^{h2}_{hhhh_-}(u, v)$, we take $L^\infty_{x,t}$, $L^2_{x,t}$, $L^4_{x,t}$ and $L^4_{x,t}$ norms to $v_{\lambda_0}$, $u_{-\lambda_0+I_{j_2}}$, $u_{\lambda_0-I_{j_1}}$ and $u_{\lambda_0-I_{j_3}}$, respectively, then
\begin{align*}
\mathscr{L}^{h2}_{hhhh_-}(u, v)&\lesssim \sum_{1\ll j_3\approx j_1\approx j_2}\langle \lambda_0\rangle^{5/4}\|\overline{v_{\lambda_0}}\|_{L^\infty_{x,t}} \|u_{-\lambda_0+I_{j_2}}\|_{L^2_{x,t}}\| u_{\lambda_0-I_{j_1}}\|_{L^4_{x,t}}\|u_{\lambda_0-I_{j_3}}\|_{L^4_{x,t}}\nonumber\\
&\lesssim \sum_{1\ll j_3\approx j_1\approx j_2}\langle \lambda_0\rangle^{5/4}\langle \lambda_0\rangle^{-1/2}(2^{j_1})^{-1/2}(2^{j_3})^{-1/2}\|v_{\lambda_0}\|_{V^2_A} \|u_{-\lambda_0+I_{j_2}}\|_{V^2_A}\nonumber\\
&\quad\quad \times T^{\varepsilon/2}\langle \lambda_0\rangle^{-3/4}(2^{j_1})^{1/4+\varepsilon}(2^{j_3})^{1/4+\varepsilon}\| u\|^2_{X^{1/4}_{\infty,A}},
\end{align*}
which is the same as the right hand side of the second inequality in \eqref{hhhh-2}.

If $u_{\lambda_0-I_{j_3}}$ has the highest dispersion modulation, we don't need to divide $\mathscr{L}^h_{hhhh_-}(u, v)$. By H\"older's inequality, we obtain that
\begin{align*}
\mathscr{L}^h_{hhhh_-}(u, v)&\lesssim \sum_{1\ll j_3\leq j_1\approx j_2}\langle \lambda_0\rangle^{5/4}\|\overline{v_{\lambda_0}}\|_{L^\infty_{x,t}} \|u_{\lambda_0-I_{j_3}}\|_{L^2_{x,t}}\| u_{\lambda_0-I_{j_1}}\|_{L^4_{x,t}}\|u_{-\lambda_0+I_{j_2}}\|_{L^4_{x,t}}\nonumber\\
&\lesssim \sum_{1\ll j_3\leq j_1\approx j_2}\langle \lambda_0\rangle^{5/4}\langle \lambda_0\rangle^{-1/2}(2^{j_1})^{-1/2}(2^{j_3})^{-1/2}\|v_{\lambda_0}\|_{V^2_A} \|u_{\lambda_0-I_{j_3}}\|_{V^2_A}\nonumber\\
&\quad\quad \times T^{\varepsilon/2}\langle \lambda_0\rangle^{-3/4}(2^{j_1})^{1/4+\varepsilon}(2^{j_2})^{1/4+\varepsilon}\| u\|^2_{X^{1/4}_{\infty,A}}\nonumber\\
&\lesssim T^{\varepsilon/2}\langle \lambda_0\rangle^{-1/4} \bigg(\sum_{j_1}(2^{j_1})^{2\varepsilon}\cdot j_1\bigg)\|v_{\lambda_0}\|_{V^2_A}\| u\|^3_{X^{1/4}_{\infty,A}}\nonumber\\
&\lesssim T^{\varepsilon/2}\|v_{\lambda_0}\|_{V^2_A}\| u\|^3_{X^{1/4}_{\infty,A}},
\end{align*}
where the last inequality is by using $2^{j_1}\lesssim \langle \lambda_0\rangle$, $j_1\lesssim (2^{j_1})^\varepsilon$, and taking $\varepsilon\leq 1/12$.

{\bf Case 2}: $\lambda_3\in h$ and $\lambda_1\in l$. In view of \eqref{FCCl}, we see that $\lambda_2\in[-3\lambda_0/4-l, \lambda_0/4-l]$, i.e., $\lambda_2\in l$ or $\lambda_2\in l_-$. We denote by $(\lambda_k)\in hhll$ that all $\lambda_0,...,\lambda_3$ satisfy the conditions \eqref{FCC}, \eqref{step1a} and
\begin{align*}
\lambda_3 \in h, \quad   \lambda_1, \lambda_2 \in l.
\end{align*}
Taking the similar notations to $hhll_{-}$,  then we divide Case 2 into two subcases.

Case $hhll$. We decompose $\lambda_1, \lambda_2, \lambda_3$ by:
\begin{align}
\lambda_k \in [0, 3\lambda_0/4]= \bigcup_{j_k\geq 0}  I_{j_k}, \ k=1,2;\ \   \lambda_3 \in [3\lambda_0/4, \lambda_0]= \bigcup_{j_3\geq 0}  \lambda_0-I_{j_3},\quad  j_1, j_2, j_3 \leq \log_2{\lambda_0}.    \nonumber
\end{align}
In view of the condition (FCC) $\lambda_1+\lambda_2+\lambda_3\approx\lambda_0$, we see that $2^{j_3}\approx 2^{j_1}+2^{j_2}$. Moreover, we can get $j_1\geq j_2$ from $\lambda_1\geq \lambda_2$. Therefore, we know that $j_3\approx j_1\geq j_2$. It means that we need to estimate
\begin{align*}
 \mathscr{L}_{hhll}(u, v):=\sum_{0\leq j_2\leq j_1\approx j_3\leq\log_2{\lambda_0}} \langle \lambda_0\rangle^{1/4}\int_{[0,T]\times\mathbb{R}} |\overline{v_{\lambda_0}} u_{I_{j_1}}u_{I_{j_2}}\partial_{x}u_{\lambda_0-I_{j_3}} |\, dxdt.
\end{align*}
In view of (DMCC) \eqref{DMCC}, we have the highest dispersion modulation
\begin{align*}
\max_{0\leq k \leq 3} |\xi^3_k - \tau_k|  \gtrsim \langle\lambda_0\rangle^2\cdot 2^{j_3}.
\end{align*}

If $v_{\lambda_0}$ has the highest dispersion modulation, by H\"older's inequality we have
\begin{align*}
 \mathscr{L}_{hhll}(u, v)\lesssim \sum_{0\leq j_2\leq j_1\approx j_3\leq\log_2{\lambda_0}} \langle \lambda_0\rangle^{5/4}\|\overline{v_{\lambda_0}}\|_{L^2_tL^{\infty}_x} \|u_{I_{j_1}}\|_{L^\infty_tL^{2}_x}\| u_{I_{j_2}}\|_{L^4_{x,t}}\|u_{\lambda_0-I_{j_3}}\|_{L^4_{x,t}}.
\end{align*}
Using $\|v_{\lambda_0}\|_{L^\infty_x}\lesssim \|v_{\lambda_0}\|_{L^2_x}$, $V^2_A\subset L^\infty_t L^2_x$, the dispersion modulation decay \eqref{dispersiondecay}, the $L^4$ estimate \eqref{lebesgue4a} and Lemma \ref{V2toX}, we have
\begin{align}
\mathscr{L}_{hhll}(u, v)&\lesssim \sum_{0\leq j_2\leq j_1\approx j_3\leq\log_2{\lambda_0}} \langle \lambda_0\rangle^{5/4}\langle \lambda_0\rangle^{-1}(2^{j_3})^{-1/2}\|v_{\lambda_0}\|_{V^2_A} \|u_{I_{j_1}}\|_{V^2_A}\nonumber\\
&\quad\quad \times T^{\varepsilon/2}(2^{j_2})^{-1/8+\varepsilon}(2^{j_3})^{1/4+\varepsilon}\langle \lambda_0\rangle^{-3/8}\| u\|^2_{X^{1/4}_{\infty,A}}\nonumber\\
&\lesssim T^{\varepsilon/2}\sum_{0\leq j_2\leq j_1\leq\log_2{\lambda_0}} \langle \lambda_0\rangle^{-1/8}(2^{j_1})^{\varepsilon}(2^{j_2})^{-1/8+\varepsilon}\|v_{\lambda_0}\|_{V^2_A}\| u\|^3_{X^{1/4}_{\infty,A}}\nonumber\\
&\lesssim T^{\varepsilon/2} \langle \lambda_0\rangle^{-1/8+\varepsilon}\|v_{\lambda_0}\|_{V^2_A}\| u\|^3_{X^{1/4}_{\infty,A}}\nonumber\\
&\lesssim T^{\varepsilon/2} \|v^{(\lambda_0)}\|_{V^2_A}\| u\|^3_{X^{1/4}_{\infty,A}},\label{step1c}
\end{align}
where the last but one inequality is obtained by summarizing over $j_2$, $j_1$ and taking $\varepsilon<1/8$.

If $u_{I_{j_1}}$ has the highest dispersion modulation, we have
\begin{align}
 \mathscr{L}_{hhll}(u, v)&\lesssim \sum_{0\leq j_2\leq j_1\approx j_3\leq\log_2{\lambda_0}} \langle \lambda_0\rangle^{5/4}\|\overline{v_{\lambda_0}}\|_{L^\infty_{x,t}} \|u_{I_{j_1}}\|_{L^2_{x,t}}\| u_{I_{j_2}}\|_{L^4_{x,t}}\|u_{\lambda_0-I_{j_3}}\|_{L^4_{x,t}}\nonumber\\
&\lesssim \sum_{0\leq j_2\leq j_1\approx j_3\leq\log_2{\lambda_0}} \langle \lambda_0\rangle^{5/4}\|v_{\lambda_0}\|_{V^2_A} \langle \lambda_0\rangle^{-1}(2^{j_3})^{-1/2}\|u_{I_{j_1}}\|_{V^2_A}\nonumber\\
&\quad\quad \times T^{\varepsilon/2}(2^{j_2})^{-1/8+\varepsilon}(2^{j_3})^{1/4+\varepsilon}\langle \lambda_0\rangle^{-3/8}\| u\|^2_{X^{1/4}_{\infty,A}},\nonumber
\end{align}
which is the same as the right hand side of the first inequality in \eqref{step1c}.

If $u_{I_{j_2}}$ has the highest dispersion modulation, we take $L^\infty_{x,t}$, $L^2_{x,t}$, $L^4_{x,t}$ and $L^4_{x,t}$ norms to $v_{\lambda_0}$, $u_{I_{j_2}}$, $u_{I_{j_1}}$ and $u_{\lambda_0-I_{j_3}}$, respectively. Then applying the dispersion modulation decay \eqref{dispersiondecay} to $u_{I_{j_2}}$, we have
\begin{align}
 \mathscr{L}_{hhll}(u, v)&\lesssim \sum_{0\leq j_2\leq j_1\approx j_3\leq\log_2{\lambda_0}} \langle \lambda_0\rangle^{5/4}\|\overline{v_{\lambda_0}}\|_{L^\infty_{x,t}} \|u_{I_{j_2}}\|_{L^2_{x,t}}\| u_{I_{j_1}}\|_{L^4_{x,t}}\|u_{\lambda_0-I_{j_3}}\|_{L^4_{x,t}}\nonumber\\
&\lesssim \sum_{0\leq j_2\leq j_1\approx j_3\leq\log_2{\lambda_0}} \langle \lambda_0\rangle^{5/4}\|v_{\lambda_0}\|_{V^2_A} \langle \lambda_0\rangle^{-1}(2^{j_3})^{-1/2}\|u_{I_{j_2}}\|_{V^2_A}\nonumber\\
&\quad\quad \times T^{\varepsilon/2}(2^{j_1})^{-1/8+\varepsilon}(2^{j_3})^{1/4+\varepsilon}\langle \lambda_0\rangle^{-3/8}\| u\|^2_{X^{1/4}_{\infty,A}}\nonumber\\
&\lesssim T^{\varepsilon/2}\langle \lambda_0\rangle^{-1/8}\sum_{0\leq j_2\leq j_1\leq\log_2{\lambda_0}}(2^{j_1})^{-3/8+2\varepsilon}(2^{j_2})^{1/4}\|v_{\lambda_0}\|_{V^2_A}\| u\|^3_{X^{1/4}_{\infty,A}}.\nonumber
\end{align}
Making the summation on $j_2$, $j_1$ in order, and taking $\varepsilon<1/16$, we can obtain the desired estimate.

If $u_{\lambda_0-I_{j_3}}$ has the highest dispersion modulation, by H\"older's inequality we have
\begin{align}
 \mathscr{L}_{hhll}(u, v)\lesssim \sum_{0\leq j_2\leq j_1\approx j_3\leq\log_2{\lambda_0}} \langle \lambda_0\rangle^{5/4}\|\overline{v_{\lambda_0}}\|_{L^\infty_{x,t}} \|u_{\lambda_0-I_{j_3}}\|_{L^2_{x,t}}\| u_{I_{j_1}}\|_{L^4_{x,t}}\| u_{I_{j_2}}\|_{L^4_{x,t}}.\nonumber
\end{align}
Using $\|v_{\lambda_0}\|_{L^\infty_x}\lesssim \|v_{\lambda_0}\|_{L^2_x}$, $V^2_A\subset L^\infty_t L^2_x$, the dispersion modulation decay \eqref{dispersiondecay}, the $L^4$ estimate \eqref{lebesgue4a} and Lemma \ref{V2toX}, we have
\begin{align}\label{step1g}
\mathscr{L}_{hhll}(u, v)&\lesssim \sum_{0\leq j_2\leq j_1\approx j_3\leq\log_2{\lambda_0}} \langle \lambda_0\rangle^{5/4}\|v_{\lambda_0}\|_{V^2_A}\langle \lambda_0\rangle^{-1}(2^{j_3})^{-1/2} \|u_{\lambda_0-I_{j_3}}\|_{V^2_A}\nonumber\\
&\quad\quad \times T^{\varepsilon/2} (2^{j_1})^{-1/8+\varepsilon}(2^{j_2})^{-1/8+\varepsilon}\| u\|^2_{X^{1/4}_{\infty,A}}\nonumber\\
&\lesssim T^{\varepsilon/2} \|v_{\lambda_0}\|_{V^2_A}\sum_{0\leq j_2\leq j_1\leq\log_2{\lambda_0}}  (2^{j_1})^{-1/8+\varepsilon} (2^{j_2})^{-1/8+\varepsilon} \| u\|^3_{X^{1/4}_{\infty,A}}\nonumber\\
&\lesssim T^{\varepsilon/2} \|v^{(\lambda_0)}\|_{V^2_A}\| u\|^3_{X^{1/4}_{\infty,A}},
\end{align}
where the last inequality is by taking $\varepsilon< 1/8$.

Case $hhll_{-}$. We decompose $\lambda_1, \lambda_2, \lambda_3$ by:
\begin{align}
\lambda_1\in [0, 3\lambda_0/4]= \bigcup_{j_1\geq 0}  I_{j_1};\ \ \lambda_2\in [-3\lambda_0/4,0]= \bigcup_{j_2\geq 0} -I_{j_2};\ \ \lambda_3 \in [3\lambda_0/4, \lambda_0]= \bigcup_{j_3\geq 0}  \lambda_0-I_{j_3}.    \nonumber
\end{align}
In view of the condition (FCC) $\lambda_1+\lambda_2+\lambda_3\approx\lambda_0$, we see that $2^{j_1}\approx 2^{j_2}+2^{j_3}$.  Thus, we know that $j_1\approx j_2\vee j_3$. If $j_1\approx j_3\geq j_2$ or $j_1\approx j_2\approx j_3$ , it is the same as Case $hhll$ to get the conclusion. So we only need to consider $j_1\approx j_2\gg j_3$,
which means that we need to estimate
\begin{align*}
 \mathscr{L}_{hhll_-}(u, v):=\sum_{j_3\ll j_2\approx j_1} \langle \lambda_0\rangle^{1/4}\int_{[0,T]\times\mathbb{R}} |\overline{v_{\lambda_0}} u_{I_{j_1}}u_{-I_{j_2}}\partial_{x}u_{\lambda_0-I_{j_3}} |\, dxdt.
\end{align*}
In view of (DMCC) \eqref{DMCC}, we have the highest dispersion modulation satisfying
\begin{align*}
\max_{0\leq k \leq 3} |\xi^3_k - \tau_k|  \gtrsim \langle\lambda_0\rangle^2\cdot 2^{j_3}.
\end{align*}

If $v_{\lambda_0}$ has the highest dispersion modulation, we can easily see that $|\lambda_0-2^{j_3}+2^{j_2}|\gtrsim \langle\lambda_0\rangle$ and $|\lambda_0-2^{j_3}-2^{j_2}|\approx|\lambda_0-2^{j_1}|\gtrsim \langle\lambda_0\rangle$, then we shall use the bilinear estimate to $u_{-I_{j_2}}u_{\lambda_0-I_{j_3}}$,
\begin{align}
 \mathscr{L}_{hhll_-}(u, v)\lesssim \sum_{j_3\ll j_2\approx j_1} \langle \lambda_0\rangle^{5/4}\|\overline{v_{\lambda_0}}\|_{L^2_tL^{\infty}_x} \|u_{I_{j_1}}\|_{L^\infty_tL^{2}_x}\| u_{-I_{j_2}}u_{\lambda_0-I_{j_3}}\|_{L^2_{x,t}}.
\end{align}
Using $\|v_{\lambda_0}\|_{L^\infty_x}\lesssim \|v_{\lambda_0}\|_{L^2_x}$, $V^2_A\subset L^\infty_t L^2_x$, the dispersion modulation decay \eqref{dispersiondecay}, the bilinear estimate \eqref{bilinear3} and Lemma \ref{V2toX}, we have
\begin{align}\label{step1b}
\mathscr{L}_{hhll_-}(u, v)&\lesssim \sum_{j_3\ll j_2\approx j_1} \langle \lambda_0\rangle^{5/4}\langle \lambda_0\rangle^{-1}(2^{j_3})^{-1/2}\|v_{\lambda_0}\|_{V^2_A} \|u_{I_{j_1}}\|_{V^2_A}\nonumber\\
&\quad\quad \times T^{\varepsilon/4}\langle \lambda_0\rangle^{-1+2\varepsilon}\| u_{-I_{j_2}}\|_{V^2_A}\|u_{\lambda_0-I_{j_3}}\|_{V^2_A}\nonumber\\
&\lesssim T^{\varepsilon/4}\sum_{j_3\ll j_2\approx j_1} \langle \lambda_0\rangle^{-3/4+2\varepsilon}(2^{j_3})^{-1/2}\|v_{\lambda_0}\|_{V^2_A} (2^{j_1})^{1/4} (2^{j_2})^{1/4} (2^{j_3})^{1/2} \langle \lambda_0\rangle^{-1/4}\| u\|^3_{X^{1/4}_{\infty,A}}\nonumber\\
&\lesssim T^{\varepsilon/4} \langle \lambda_0\rangle^{-1+2\varepsilon}\bigg(\sum_{j_1}(2^{j_1})^{1/2}\cdot j_1\bigg) \|v_{\lambda_0}\|_{V^2_A}\| u\|^3_{X^{1/4}_{\infty,A}}\nonumber\\
&\lesssim T^{\varepsilon/4} \|v^{(\lambda_0)}\|_{V^2_A}\| u\|^3_{X^{1/4}_{\infty,A}},
\end{align}
where the last but one inequality is obtained by taking $\varepsilon<1/4$.

If $u_{I_{j_1}}$ has the highest dispersion modulation, we have
\begin{align}
 \mathscr{L}_{hhll_-}(u, v)&\lesssim \sum_{j_3\ll j_2\approx j_1} \langle \lambda_0\rangle^{5/4}\|\overline{v_{\lambda_0}}\|_{L^\infty_{x,t}} \|u_{I_{j_1}}\|_{L^2_{x,t}}\| u_{-I_{j_2}}u_{\lambda_0-I_{j_3}}\|_{L^2_{x,t}}\nonumber\\
&\lesssim \sum_{j_3\ll j_2\approx j_1} \langle \lambda_0\rangle^{5/4}\|v_{\lambda_0}\|_{V^2_A} \langle \lambda_0\rangle^{-1}(2^{j_3})^{-1/2}\|u_{I_{j_1}}\|_{V^2_A}\nonumber\\
&\quad\quad \times T^{\varepsilon/4}\langle \lambda_0\rangle^{-1+2\varepsilon}\| u_{-I_{j_2}}\|_{V^2_A}\|u_{\lambda_0-I_{j_3}}\|_{V^2_A},\nonumber
\end{align}
which is the same as the right hand side of the first inequality in \eqref{step1b}.

If $u_{-I_{j_2}}$ has the highest dispersion modulation, noticing that $j_1\approx j_2$, we can take $L^\infty_{x,t}$, $L^2_{x,t}$ and $L^2_{x,t}$ norms to $v_{\lambda_0}$, $u_{-I_{j_2}}$ and $u_{I_{j_1}}u_{\lambda_0-I_{j_3}}$, respectively. Then we can repeat the above proof to obtain the desired estimates.

If $u_{\lambda_0-I_{j_3}}$ has the highest dispersion modulation, comparing with Case $hhll$, the difference is the summation in \eqref{step1g} (taking $\varepsilon< 1/8$)
\begin{align}
\sum_{j_3\ll j_2\approx j_1}  (2^{j_1})^{-1/8+\varepsilon} (2^{j_2})^{-1/8+\varepsilon}\lesssim \sum_{0\leq j_1\leq\log_2{\lambda_0}}  (2^{j_1})^{-1/4+2\varepsilon} \times j_1\lesssim 1.
\end{align}

{\bf Case 3}: $\lambda_3\in h$ and $\lambda_1\in l_-$. It is easy to see that $\lambda_2\in l_-$. We decompose $\lambda_1, \lambda_2, \lambda_3$ by:
\begin{align}
\lambda_k\in [-c\lambda_0,0]= \bigcup_{j_k\geq 0} -I_{j_k}, k=1,2;\ \ \lambda_3 \in [c\lambda_0, \lambda_0]= \bigcup_{j_3\geq 0}  \lambda_0-I_{j_3}.    \nonumber
\end{align}
In view of the condition (FCC) $\lambda_1+\lambda_2+\lambda_3\approx\lambda_0$, we see that $2^{j_1}+2^{j_2}+2^{j_3}\approx 0$. It means that $0\leq j_1,j_2,j_3\lesssim 1$. Then using the dispersion modulation decay \eqref{dispersiondecay}, the $L^4$ estimate \eqref{lebesgue4a} and Lemma \ref{V2toX}, and noticing that the summation about $j_1,\ j_2,\ j_3$ is finite, we can get the result and the details are omitted.

{\bf Case 4}: $\lambda_3\in l$. This case is easy to estimate because the derivative locates in the low frequency, $\lambda_1,\lambda_2\in \{l,l_-\}$ and the highest dispersion modulation satisfies
\begin{align*}
\max_{0\leq k \leq 3} |\xi^3_k - \tau_k|  \gtrsim \langle\lambda_0\rangle^3.
\end{align*}
We take Case $hlll$ ($\lambda_1,\lambda_2,\lambda_3\in l$) as an example. When $v_{\lambda_0}$ attains the highest dispersion modulation, using a similar way as above, we have
\begin{align}
 \mathscr{L}_{hlll}(u, v)&\lesssim \sum_{j_1,j_2,j_3} \langle \lambda_0\rangle^{1/4} 2^{j_3}\|\overline{v_{\lambda_0}}\|_{L^2_tL^{\infty}_x} \|u_{I_{j_1}}\|_{L^\infty_tL^{2}_x}\| u_{I_{j_2}}\|_{L^4_{x,t}}\|u_{I_{j_3}}\|_{L^4_{x,t}}\nonumber\\
&\lesssim \sum_{j_1,j_2,j_3} \langle \lambda_0\rangle^{1/4} 2^{j_3}\langle \lambda_0\rangle^{-3/2}\|v_{\lambda_0}\|_{V^2_A}(2^{j_1})^{1/4}T^{\varepsilon/2}(2^{j_2})^{-1/8+\varepsilon} (2^{j_3})^{-1/8+\varepsilon}\| u\|^3_{X^{1/4}_{\infty,A}}\nonumber\\
&\lesssim T^{\varepsilon/2}\langle \lambda_0\rangle^{-5/4}\sum_{j_1,j_2,j_3}(2^{j_1})^{1/4}(2^{j_2})^{-1/8+\varepsilon}(2^{j_3})^{7/8+\varepsilon}\|v_{\lambda_0}\|_{V^2_A}\| u\|^3_{X^{1/4}_{\infty,A}}\nonumber\\
&\lesssim T^{\varepsilon/2}\|v_{\lambda_0}\|_{V^2_A}\| u\|^3_{X^{1/4}_{\infty,A}}.\nonumber
\end{align}
When $u_{I_{j_1}}$, $u_{I_{j_2}}$, or $u_{I_{j_3}}$ attains the highest dispersion modulation, we can use an analogous way to get the result. In fact, we just need to take $L^\infty_{x,t}$ norm to $v_{\lambda_0}$, $L^2_{x,t}$ norm to the item which has the highest dispersion modulation, and $L^4_{x,t}$ norm to the other two items.

{\bf Step 2.} We consider the case that $\lambda_0$ is the secondly maximal integer in $\lambda_0, \cdots, \lambda_3$. By the symmetry, we can assume $\lambda_1\geq\lambda_2$. Then $\lambda_0, \cdots, \lambda_3$ have the following three orders:
\begin{align*}
{\rm Order\,1:}\quad\quad\lambda_3\geq \lambda_0 \geq \lambda_1 \geq \lambda_2;\\
{\rm Order\,2:}\quad\quad\lambda_1\geq \lambda_0 \geq \lambda_3 \geq \lambda_2;\\
{\rm Order\,3:}\quad\quad\lambda_1\geq \lambda_0 \geq \lambda_2 \geq \lambda_3.
\end{align*}
Considering the derivative is located in $u_{\lambda_3}$, we take the Order 1 for example in the following proof (the other orders are similar). We divide the proof into three cases $|\lambda_0|\lesssim 1$, $\lambda_0\ll 0$ and $\lambda_0 \gg 0$.

{\bf Case 1:} $|\lambda_0|\lesssim 1$. We decompose $\lambda_1, \lambda_2, \lambda_3$ by:
\begin{align}
\lambda_k\in (-\infty,\lambda_0]= \bigcup_{j_k\geq -1} -I_{j_k}, k=1,2;\ \ \lambda_3 \in [\lambda_0, +\infty)= \bigcup_{j_3\geq -1}  I_{j_3},\ \  I_{-1}=[-|\lambda_0|,0).    \nonumber
\end{align}
In view of $\lambda_0\approx\lambda_1+\lambda_2+\lambda_3$ and $\lambda_1\geq\lambda_2$, we have $j_3\approx j_2\geq j_1\geq -1$. By DMCC \eqref{DMCC} the highest dispersion modulation satisfies
\begin{align}
\max_{0\leq k \leq 3} |\xi^3_k - \tau_k|  \gtrsim 2^{j_1}\cdot 2^{j_2}\cdot 2^{j_3}.
\end{align}

If $v_{\lambda_0}$ gains the highest dispersion modulation, we have
\begin{align}\label{step2a}
 &\sum_{j_3\approx j_2\geq j_1\geq -1} \langle \lambda_0\rangle^{1/4}\int_{[0,T]\times\mathbb{R}} |\overline{v_{\lambda_0}} u_{-I_{j_1}}u_{-I_{j_2}}\partial_{x}u_{I_{j_3}}| \, dxdt\nonumber\\
 \lesssim & \sum_{j_3\approx j_2\geq j_1\geq -1} 2^{j_3}\|\overline{v_{\lambda_0}}\|_{L^2_tL^{\infty}_x} \|u_{-I_{j_1}}\|_{L^\infty_tL^{2}_x}\| u_{-I_{j_2}}\|_{L^4_{x,t}}\|u_{I_{j_3}}\|_{L^4_{x,t}}\nonumber\\
\lesssim & \sum_{j_3\approx j_2\geq j_1\geq -1} 2^{j_3}(2^{j_1})^{-1/2}(2^{j_2})^{-1/2}(2^{j_3})^{-1/2}\|v_{\lambda_0}\|_{V^2_A} \|u_{-I_{j_1}}\|_{V^2_A}\nonumber\\
&\quad \times T^{\varepsilon/4}(2^{j_2})^{-1/8+\varepsilon} T^{\varepsilon/4}(2^{j_3})^{-1/8+\varepsilon}\| u\|^2_{X^{1/4}_{\infty,A}}\nonumber\\
\lesssim & T^{\varepsilon/2}\sum_{j_3\geq j_1\geq -1} (2^{j_3})^{-1/4+2\varepsilon}(2^{j_1})^{-1/4}\|v_{\lambda_0}\|_{V^2_A}\|u\|^3_{X^{1/4}_{\infty,A}}\nonumber\\
\lesssim &T^{\varepsilon/2} \|v^{(\lambda_0)}\|_{V^2_A}\|u\|^3_{X^{1/4}_{\infty,A}}.
\end{align}

If $u_{-I_{j_1}}$ has the highest dispersion modulation, we take $L^\infty_{x,t}$, $L^2_{x,t}$, $L^4_{x,t}$  and $L^4_{x,t}$ norms to $v_{\lambda_0}$, $u_{-I_{j_1}}$, $u_{-I_{j_2}}$ and $u_{I_{j_3}}$, respectively. Then applying the dispersion modulation decay \eqref{dispersiondecay}  and the $L^4$ estimate Lemma \ref{L4}, we can get the desired conclusion.

If $u_{I_{j_3}}$ has the highest dispersion modulation, we divide the left hand side of \eqref{trilinear4} into two terms.
\begin{align}
&\sum_{j_3\approx j_2\geq j_1\geq -1} \langle \lambda_0\rangle^{1/4}\int_{[0,T]\times\mathbb{R}} |\overline{v_{\lambda_0}} u_{-I_{j_1}}u_{-I_{j_2}}\partial_{x}u_{I_{j_3}}| \, dxdt\nonumber\\
\leq &\bigg(\sum_{j_3\approx j_2\approx j_1\geq -1}+\sum_{j_3\approx j_2\gg j_1\geq -1}\bigg)\langle\lambda_0\rangle^{1/4}\int_{[0,T]\times\mathbb{R}} |\overline{v_{\lambda_0}} u_{-I_{j_1}}u_{-I_{j_2}}\partial_{x}u_{I_{j_3}}| \, dxdt\nonumber\\
:=&I_1(u,v)+I_2(u,v).
\end{align}
For $I_1(u,v)$, $L^4$ estimate \eqref{lebesgue4a} is enough.
\begin{align}
I_1(u,v)&\lesssim \sum_{j_3\approx j_2\approx j_1\geq -1} 2^{j_3}\|\overline{v_{\lambda_0}}\|_{L^\infty_{x,t}} \|u_{I_{j_3}}\|_{L^2_{x,t}}\| u_{-I_{j_1}}\|_{L^4_{x,t}}\|u_{-I_{j_2}}\|_{L^4_{x,t}}\nonumber\\
&\lesssim \sum_{j_3\approx j_2\approx j_1\geq -1}2^{j_3}\|v_{\lambda_0}\|_{V^2_A} (2^{j_1})^{-1/2}(2^{j_2})^{-1/2}(2^{j_3})^{-1/2}(2^{j_3})^{1/4}\|u\|_{X^{1/4}_{\infty,A}}\nonumber\\
&\quad\times T^{\varepsilon/4}(2^{j_1})^{-1/8+\varepsilon} T^{\varepsilon/4}(2^{j_2})^{-1/8+\varepsilon}\| u\|^2_{X^{1/4}_{\infty,A}}\nonumber\\
&\lesssim T^{\varepsilon/2}\sum_{j_3\geq -1}(2^{j_3})^{-1/2+2\varepsilon}\|v_{\lambda_0}\|_{V^2_A}\|u\|^3_{X^{1/4}_{\infty,A}}\nonumber\\
&\lesssim T^{\varepsilon/2}\|v^{(\lambda_0)}\|_{V^2_A}\|u\|^3_{X^{1/4}_{\infty,A}}.
\end{align}
For $I_2(u,v)$, we need to use the bilinear estimate \eqref{bilinear3}.
\begin{align}
I_2(u,v)&\lesssim \sum_{j_3\approx j_2\gg j_1\geq -1} 2^{j_3}\|\overline{v_{\lambda_0}}\|_{L^\infty_{x,t}} \|u_{I_{j_3}}\|_{L^2_{x,t}}\| u_{-I_{j_1}}u_{-I_{j_2}}\|_{L^2_{x,t}}\nonumber\\
&\lesssim \sum_{j_3\approx j_2\gg j_1\geq -1}2^{j_3}\|v_{\lambda_0}\|_{V^2_A} (2^{j_1})^{-1/2}(2^{j_2})^{-1/2}(2^{j_3})^{-1/2}(2^{j_3})^{1/4}\|u\|_{X^{1/4}_{\infty,A}}\nonumber\\
&\quad\times T^{\varepsilon/4}(2^{j_2})^{-1+2\varepsilon} (2^{j_1})^{1/4}(2^{j_2})^{1/4}\| u\|^2_{X^{1/4}_{\infty,A}}\nonumber\\
&\lesssim T^{\varepsilon/4}\sum_{j_3\gg j_1\geq -1}(2^{j_3})^{-1/2+2\varepsilon}(2^{j_1})^{-1/4}\|v_{\lambda_0}\|_{V^2_A}\|u\|^3_{X^{1/4}_{\infty,A}}\nonumber\\
&\lesssim T^{\varepsilon/4}\|v^{(\lambda_0)}\|_{V^2_A}\|u\|^3_{X^{1/4}_{\infty,A}}.
\end{align}

If $u_{-I_{j_2}}$ has the highest dispersion modulation,  we can get the desired estimate by exchanging the positions of $u_{I_{j_3}}$ and $u_{-I_{j_2}}$ in the above discussion(noticing that $j_2\approx j_3$).

{\bf Case 2:} $\lambda_0\ll 0$. We decompose $\lambda_1$ and $\lambda_2$ by:
\begin{align}
\lambda_k\in (-\infty,\lambda_0]= \bigcup_{j_k\geq 0} \lambda_0-I_{j_k},\quad k=1,2. \nonumber
\end{align}
From the following frequency constraint condition
\begin{align}
\lambda_0=\lambda_1+\lambda_2+\lambda_3+l,\quad |l|\leq 10,
\end{align}
we can decompose $\lambda_3$ as follows.
\begin{align}
\lambda_3\in [-\lambda_0-l,+\infty]= \bigcup_{j_3\geq -1} -\lambda_0+I_{j_3},\quad I_{-1}=[-|l|,0). \nonumber
\end{align}
In view of $\lambda_0\approx\lambda_1+\lambda_2+\lambda_3$ and $\lambda_1\geq\lambda_2$, we have $j_3\approx j_2\geq j_1$. By DMCC \eqref{DMCC}, we can see that the highest dispersion modulation satisfies
\begin{align}
\max_{0\leq k \leq 3} |\xi^3_k - \tau_k|  \gtrsim 2^{j_1}\cdot 2^{j_2}\cdot (\langle\lambda_0\rangle+2^{j_3}).
\end{align}

If the highest dispersion modulation is located in $v_{\lambda_0}$, from the dispersion modulation decay \eqref{dispersiondecay}, $L^4$ estimate \eqref{lebesgue4a} and Lemma \ref{V2toX}, we have
\begin{align}\label{step2b}
&\sum_{j_3\approx j_2\geq j_1} \langle \lambda_0\rangle^{1/4}\int_{[0,T]\times\mathbb{R}} |\overline{v_{\lambda_0}} u_{\lambda_0-I_{j_1}}u_{\lambda_0-I_{j_2}}\partial_{x}u_{-\lambda_0+I_{j_3}}| \, dxdt\nonumber\\
\lesssim &\sum_{j_3\approx j_2\geq j_1}\langle \lambda_0\rangle^{1/4}(\langle\lambda_0\rangle+ 2^{j_3})\|\overline{v_{\lambda_0}}\|_{L^2_tL^{\infty}_x} \|u_{\lambda_0-I_{j_1}}\|_{L^\infty_tL^{2}_x}\| u_{\lambda_0-I_{j_2}}\|_{L^4_{x,t}}\|u_{-\lambda_0+I_{j_3}}\|_{L^4_{x,t}}\nonumber\\
\lesssim &\sum_{j_3\approx j_2\geq j_1}\langle \lambda_0\rangle^{1/4}(\langle\lambda_0\rangle+ 2^{j_3})(2^{j_2})^{-1/2} (\langle\lambda_0\rangle+2^{j_3})^{-1/2}\|v_{\lambda_0}\|_{V^2_A} (\langle\lambda_0\rangle+2^{j_1})^{-1/4}\nonumber\\
&\quad\times T^{\varepsilon/4}(2^{j_2})^{1/4+\varepsilon}(\langle\lambda_0\rangle+2^{j_2})^{-3/8} T^{\varepsilon/4}(2^{j_3})^{1/4+\varepsilon}(\langle\lambda_0\rangle+2^{j_3})^{-3/8}\| u\|^3_{X^{1/4}_{\infty,A}}\nonumber\\
\lesssim &T^{\varepsilon/2}\sum_{j_3\geq 0}(\langle\lambda_0\rangle+2^{j_3})^{-1/4}(2^{j_3})^{2\varepsilon}\sum_{0\leq j_1\leq j_3}\langle \lambda_0\rangle^{1/4}(\langle\lambda_0\rangle+2^{j_1})^{-1/4}\|v_{\lambda_0}\|_{V^2_A}\|u\|^3_{X^{1/4}_{\infty,A}}.
\end{align}
Making the summation on $j_1$, we see that the summation is controlled by $j_3$. Then one has that for $0<\varepsilon<1/8$,
\begin{align}
\eqref{step2b}\lesssim &T^{\varepsilon/2}\bigg(\sum_{j_3\geq 0}(2^{j_3})^{-1/4+2\varepsilon}\cdot j_3\bigg)\|v_{\lambda_0}\|_{V^2_A}\|u\|^3_{X^{1/4}_{\infty,A}}\nonumber\\
\lesssim & T^{\varepsilon/2}\|v^{(\lambda_0)}\|_{V^2_A}\|u\|^3_{X^{1/4}_{\infty,A}}.\nonumber
\end{align}

If the highest dispersion modulation is located in $u_{\lambda_0-I_{j_1}}$, we take $L^\infty_{x,t}$, $L^2_{x,t}$, $L^4_{x,t}$  and $L^4_{x,t}$ norms to $v_{\lambda_0}$, $u_{\lambda_0-I_{j_1}}$, $u_{\lambda_0-I_{j_2}}$ and $u_{-\lambda_0+I_{j_3}}$, respectively. Then we can reduce the  desired estimate as the above case, so the details are omitted.

If the highest dispersion modulation is located in $u_{\lambda_0-I_{j_2}}$, we divide the left hand side of \eqref{trilinear4} into two terms.
\begin{align}
&\sum_{j_3\approx j_2\geq j_1\geq 0} \langle \lambda_0\rangle^{1/4}\int_{[0,T]\times\mathbb{R}} |\overline{v_{\lambda_0}} u_{\lambda_0-I_{j_1}}u_{\lambda_0-I_{j_2}}\partial_{x}u_{-\lambda_0+I_{j_3}}| \, dxdt\nonumber\\
\leq &\bigg(\sum_{j_3\approx j_2\approx j_1\geq 0}+\sum_{j_3\approx j_2\gg j_1\geq 0}\bigg)\langle \lambda_0\rangle^{1/4}\int_{[0,T]\times\mathbb{R}} |\overline{v_{\lambda_0}} u_{\lambda_0-I_{j_1}}u_{\lambda_0-I_{j_2}}\partial_{x}u_{-\lambda_0+I_{j_3}}| \, dxdt\nonumber\\
:=&I_1(u,v)+I_2(u,v).
\end{align}
For $I_1(u,v)$, from the dispersion modulation decay \eqref{dispersiondecay}, $L^4$ estimate \eqref{lebesgue4a} and Lemma \ref{V2toX}, we have
\begin{align}
I_1(u,v)&\lesssim \sum_{j_3\approx j_2\approx j_1\geq 0} \langle \lambda_0\rangle^{1/4}(\langle\lambda_0\rangle+ 2^{j_3})\|\overline{v_{\lambda_0}}\|_{L^\infty_{x,t}} \|u_{\lambda_0-I_{j_2}}\|_{L^2_{x,t}}\| u_{\lambda_0-I_{j_1}}\|_{L^4_{x,t}}\|u_{-\lambda_0+I_{j_3}}\|_{L^4_{x,t}}\nonumber\\
&\lesssim \sum_{j_3\approx j_2\approx j_1\geq 0}\langle \lambda_0\rangle^{1/4}(\langle\lambda_0\rangle+ 2^{j_3})\|v_{\lambda_0}\|_{V^2_A} (2^{j_1})^{-1/2}(\langle\lambda_0\rangle+2^{j_3})^{-1/2} (\langle\lambda_0\rangle+2^{j_2})^{-1/4}\nonumber\\
&\quad\times T^{\varepsilon/4}(2^{j_1})^{1/4+\varepsilon}(\langle\lambda_0\rangle+2^{j_1})^{-3/8} T^{\varepsilon/4}(2^{j_3})^{1/4+\varepsilon}(\langle\lambda_0\rangle+2^{j_3})^{-3/8}\| u\|^3_{X^{1/4}_{\infty,A}}\nonumber\\
&\lesssim T^{\varepsilon/2}\sum_{j_3\geq 0}\langle \lambda_0\rangle^{1/4}(\langle\lambda_0\rangle+2^{j_3})^{-1/2}(2^{j_3})^{2\varepsilon}\|v_{\lambda_0}\|_{V^2_A}\|u\|
^3_{X^{1/4}_{\infty,A}}.\label{step2c}
\end{align}
Noticing that
\begin{align}
\langle \lambda_0\rangle^{1/4}(\langle\lambda_0\rangle+2^{j_3})^{-1/4}\leq 1, \quad (\langle\lambda_0\rangle+2^{j_3})^{-1/4}(2^{j_3})^{2\varepsilon}\leq (2^{j_3})^{-1/4+2\varepsilon},\nonumber
\end{align}
for $0<\varepsilon<1/8$, \eqref{step2c} is dominated by
\begin{align}
&\lesssim T^{\varepsilon/2}\sum_{j_3\geq 0}(2^{j_3})^{-1/4+2\varepsilon}\|v_{\lambda_0}\|_{V^2_A}\|u\|^3_{X^{1/4}_{\infty,A}}\lesssim T^{\varepsilon/2}\|v^{(\lambda_0)}\|_{V^2_A}\|u\|^3_{X^{1/4}_{\infty,A}}.\nonumber
\end{align}
For $I_2(u,v)$, from the dispersion modulation decay \eqref{dispersiondecay}, the bilinear estimate \eqref{bilinear3} and Lemma \ref{V2toX}, we have
\begin{align}
I_2(u,v)&\lesssim \sum_{j_3\approx j_2\gg j_1\geq 0} \langle \lambda_0\rangle^{1/4}(\langle\lambda_0\rangle+ 2^{j_3})\|\overline{v_{\lambda_0}}\|_{L^\infty_{x,t}} \|u_{\lambda_0-I_{j_2}}\|_{L^2_{x,t}}\| u_{\lambda_0-I_{j_1}}u_{-\lambda_0+I_{j_3}}\|_{L^2_{x,t}}\nonumber\\
&\lesssim \sum_{j_3\approx j_2\gg j_1\geq 0}\langle \lambda_0\rangle^{1/4} (2^{j_1})^{-1/2}(2^{j_2})^{-1/2} (\langle\lambda_0\rangle+2^{j_3})^{1/2}\|v_{\lambda_0}\|_{V^2_A}\|u_{\lambda_0-I_{j_2}}\|_{V^2_A}\nonumber\\
&\quad\times T^{\varepsilon/4}(\langle\lambda_0\rangle+2^{j_3})^{-1/2+\varepsilon}(2^{j_3})^{-1/2+\varepsilon}\|u
_{\lambda_0-I_{j_1}}\|_{V^2_A} \|u_{-\lambda_0+I_{j_3}}\|_{V^2_A}\nonumber\\
&\lesssim T^{\varepsilon/4}\sum_{j_3\gg j_1\geq 0}\langle \lambda_0\rangle^{1/4}(\langle\lambda_0\rangle+2^{j_3})^{-1/2+\varepsilon}(2^{j_3})^{\varepsilon} (\langle\lambda_0\rangle+2^{j_1})^{-1/4}\|v_{\lambda_0}\|_{V^2_A}\|u\|^3_{X^{1/4}_{\infty,A}}\nonumber\\
&\lesssim T^{\varepsilon/4}\bigg(\sum_{j_3\geq 0}(2^{j_3})^{-1/2+2\varepsilon}\cdot j_3\bigg) \|v_{\lambda_0}\|_{V^2_A}\|u\|^3_{X^{1/4}_{\infty,A}}\nonumber\\
&\lesssim T^{\varepsilon/4}\|v^{(\lambda_0)}\|_{V^2_A}\|u\|^3_{X^{1/4}_{\infty,A}}.\nonumber
\end{align}

If the highest dispersion modulation is located in $u_{-\lambda_0+I_{j_3}}$, noticing that $j_2\approx j_3$ and $|-\lambda_0+I_{j_3}|\sim |\lambda_0-I_{j_2}|\sim (\langle\lambda_0\rangle+2^{j_3})$, we can get the desired estimate by exchanging the positions of $u_{-\lambda_0+I_{j_3}}$ and $u_{\lambda_0-I_{j_2}}$ in the above discussion.

{\bf Case 3:} $\lambda_0\gg 0$. From the frequency constraint condition $\lambda_0=\lambda_1+\lambda_2+\lambda_3+l,\quad |l|\leq 10$, we know that $\lambda_2$ must be less than zero. Furthermore, one can divide this case into three subcases:
$\lambda_2\in [-c\lambda_0, 0]$, $\lambda_2\in [-\lambda_0, -c\lambda_0]$ and $\lambda_2\in (-\infty,-\lambda_0]$.

{\bf Case 3.1:} $\lambda_2\in [-c\lambda_0, 0]$. From the frequency constraint condition we find that $\lambda_1\in[-c\lambda_0,0]$ or $[0,c\lambda_0]$ ($\lambda_1\in[c\lambda_0,\lambda_0]$ will never happen), and $\lambda_3$ satisfies Table \ref{table0}.
\begin{table}
\begin{center}
\begin{tabular}{|c|c|c|c|}
\hline
 ${\rm Case}$   &  $\lambda_2\in $ & $\lambda_1\in $ & $\lambda_3\in $    \\
\hline
$l_-l_- h$  & $ [-c\lambda_0,0]$  & $ [-c\lambda_0, 0]$  &  $[\lambda_0, \lambda_0+2c\lambda_0-l]$  \\
\hline
$l_-l h$  & $ [-c\lambda_0,0]$  & $ [0,c\lambda_0]$  &  $[\lambda_0, \lambda_0+c\lambda_0-l]$  \\
\hline
\end{tabular}
\end{center}
\caption{$\lambda_2\in [-c\lambda_0,0]$}\label{table0}
\end{table}

{\it Case $l_-l_- h$.} One can use the dyadic decomposition:
\begin{gather*}
\lambda_k\in [-c\lambda_0,0]= \bigcup_{j_k\geq 0}-I_{j_k},\ k=1,2; \\
\lambda_3 \in [\lambda_0, (1+2c)\lambda_0-l]= \bigcup_{j_3\geq 0} \lambda_0+I_{j_3}, \ \ j_1, j_2, j_3\leq \log_2 \langle\lambda_0\rangle+1.
\end{gather*}
From the frequency constraint condition \eqref{FCC}, we know
\begin{align}
2^{j_1}+2^{j_2}\approx 2^{j_3}\ \Rightarrow\ j_3\approx j_1\vee j_2.
\end{align}
We can easily get that the highest dispersion modulation satisfies
\begin{align}
\max_{0\leq k \leq 3} |\xi^3_k - \tau_k|  \gtrsim \langle\lambda_0\rangle^2\cdot 2^{j_3}.
\end{align}

If $v_{\lambda_0}$ attains the highest dispersion modulation, from the dispersion modulation decay \eqref{dispersiondecay}, $L^4$ estimate \eqref{lebesgue4a} and Lemma \ref{V2toX}, we have
\begin{align*}
 &\sum_{j_3\approx j_1\vee j_2} \langle \lambda_0\rangle^{1/4}\int_{[0,T]\times\mathbb{R}} |\overline{v_{\lambda_0}} u_{-I_{j_1}}u_{-I_{j_2}}\partial_{x}u_{\lambda_0+I_{j_3}}| \, dxdt\nonumber\\
 \lesssim & \sum_{j_3\approx j_1\vee j_2} \langle \lambda_0\rangle^{5/4}\|\overline{v_{\lambda_0}}\|_{L^2_tL^{\infty}_x} \| u_{-I_{j_1}}\|_{L^\infty_tL^2_x} \|u_{-I_{j_2}}\|_{L^4_{x,t}}\|u_{\lambda_0+I_{j_3}}\|_{L^4_{x,t}}\nonumber\\
\lesssim & T^{\varepsilon/2}\sum_{j_3\approx j_1\vee j_2} \langle \lambda_0\rangle^{-1/8}(2^{j_3})^{-1/4+\varepsilon}(2^{j_1})^{1/4}
(2^{j_2})^{-1/8+\varepsilon}\|v_{\lambda_0}\|_{V^2_A}\|u\|^3_{X^{1/4}_{\infty,A}}\nonumber\\
\lesssim &T^{\varepsilon/2}\langle \lambda_0\rangle^{-1/8+\varepsilon}  \|v_{\lambda_0}\|_{V^2_A}\|u\|^3_{X^{1/4}_{\infty,A}}\nonumber\\
\lesssim &T^{\varepsilon/2} \|v^{(\lambda_0)}\|_{V^2_A}\|u\|^3_{X^{1/4}_{\infty,A}},
\end{align*}
where the last but one inequality is gained by summarizing over $j_2,j_1$ and $j_3$ in order. One just note that $j_1\leq j_3\leq \log_2 \langle\lambda_0\rangle+1$ and take $\varepsilon<1/8$.

If $u_{-I_{j_1}}$ has the highest dispersion modulation,  we take $L^\infty_{x,t}$, $L^2_{x,t}$, $L^4_{x,t}$ and $L^4_{x,t}$ norms to $v_{\lambda_0}$, $u_{-I_{j_1}}$, $u_{-I_{j_2}}$ and $u_{\lambda_0+I_{j_3}}$, respectively. Then we can get the desired conclusion by the same way as above. If $u_{-I_{j_2}}$ gains the highest dispersion modulation, one can exchange the positions of $j_1$ and $j_2$ to obtain the desired estimate.

If $u_{\lambda_0+I_{j_3}}$ has the highest dispersion modulation, we have
\begin{align}\label{step2n}
 &\sum_{j_3\approx j_1\vee j_2} \langle \lambda_0\rangle^{1/4}\int_{[0,T]\times\mathbb{R}} |\overline{v_{\lambda_0}} u_{-I_{j_1}}u_{-I_{j_2}}\partial_{x}u_{\lambda_0+I_{j_3}}| \, dxdt\nonumber\\
\lesssim & \sum_{j_3\approx j_1\vee j_2} \langle \lambda_0\rangle^{5/4}\|\overline{v_{\lambda_0}}\|_{L^\infty_{x,t}}\| u_{\lambda_0+I_{j_3}}\|_{L^2_{x,t}} \|u_{-I_{j_1}}\|_{L^4_{x,t}}\|u_{-I_{j_2}}\|_{L^4_{x,t}}\nonumber\\
\lesssim & T^{\varepsilon/2}\sum_{j_3\approx j_1\vee j_2}(2^{j_1})^{-1/8+\varepsilon} (2^{j_2})^{-1/8+\varepsilon} \|v_{\lambda_0}\|_{V^2_A} \|u\|^3_{X^{1/4}_{\infty,A}}\nonumber\\
\lesssim & T^{\varepsilon/2} \|v^{(\lambda_0)}\|_{V^2_A}\|u\|^3_{X^{1/4}_{\infty,A}}.
\end{align}

{\it Case $l_-l h$.} One can use the dyadic decomposition:
\begin{gather*}
\lambda_1\in [0, c\lambda_0]= \bigcup_{j_1\geq 0}I_{j_1},\ \lambda_2\in [-c\lambda_0,0]= \bigcup_{j_2\geq 0}-I_{j_2},\\
\lambda_3 \in [\lambda_0, (1+c)\lambda_0-l]= \bigcup_{j_3\geq 0} \lambda_0+I_{j_3}, \ \ j_1, j_2, j_3\leq \log_2 \langle\lambda_0\rangle.
\end{gather*}
From the frequency constraint condition \eqref{FCC}, we get
\begin{align}
2^{j_1}+2^{j_3}\approx 2^{j_2}\ \Rightarrow\ j_2\approx j_1\vee j_3.
\end{align}
One can get that the highest dispersion modulation satisfies
\begin{align}
\max_{0\leq k \leq 3} |\xi^3_k - \tau_k|  \gtrsim \langle\lambda_0\rangle^2\cdot 2^{j_3}.
\end{align}

If $v_{\lambda_0}$ attains the highest dispersion modulation, from the dispersion modulation decay \eqref{dispersiondecay}, the bilinear estimate  and Lemma \ref{V2toX}, we have
\begin{align*}
 &\sum_{j_2\approx j_1\vee j_3} \langle \lambda_0\rangle^{1/4}\int_{[0,T]\times\mathbb{R}} |\overline{v_{\lambda_0}} u_{I_{j_1}}u_{-I_{j_2}}\partial_{x}u_{\lambda_0+I_{j_3}}| \, dxdt\nonumber\\
 \lesssim & \sum_{j_2\approx j_1\vee j_3} \langle \lambda_0\rangle^{5/4}\|\overline{v_{\lambda_0}}\|_{L^2_tL^{\infty}_x} \| u_{I_{j_1}}\|_{L^\infty_tL^2_x} \|u_{-I_{j_2}}u_{\lambda_0+I_{j_3}}\|_{L^2_{x,t}}\nonumber\\
\lesssim & T^{\varepsilon/4}\sum_{j_2\approx j_1\vee j_3} \langle \lambda_0\rangle^{1/4}(2^{j_3})^{-1/2}\|v_{\lambda_0}\|_{V^2_A}\|u_{I_{j_1}}\|_{V^2_A}\langle \lambda_0\rangle^{-1+2\varepsilon}\|u_{-I_{j_2}}\|_{V^2_A}\|u_{\lambda_0+I_{j_3}}\|_{V^2_A}\nonumber\\
\lesssim &T^{\varepsilon/4}\sum_{j_2\approx j_1\vee j_3} \langle \lambda_0\rangle^{-1+2\varepsilon}(2^{j_1})^{1/4}(2^{j_2})^{1/4}  \|v_{\lambda_0}\|_{V^2_A}\|u\|^3_{X^{1/4}_{\infty,A}}\nonumber\\
\lesssim &T^{\varepsilon/4} \|v^{(\lambda_0)}\|_{V^2_A}\|u\|^3_{X^{1/4}_{\infty,A}}.
\end{align*}

If $u_{I_{j_1}}$ has the highest dispersion modulation,  we just take $L^\infty_{x,t}$, $L^2_{x,t}$ and $L^2_{x,t}$ norms to $v_{\lambda_0}$, $u_{I_{j_1}}$ and $u_{-I_{j_2}}u_{\lambda_0+I_{j_3}}$, respectively. If $u_{-I_{j_2}}$ attains the highest dispersion modulation, one can further exchange the positions of $j_1$ and $j_2$ to obtain the desired estimate.
If $u_{\lambda_0+I_{j_3}}$ has the highest dispersion modulation, we can get the result by the same way as \eqref{step2n} in Case $l_-l_-h$.

{\bf Case 3.2:} $\lambda_2\in [-\lambda_0,-c\lambda_0]$.  We consider $\lambda_1\in [c\lambda_0,\lambda_0]$, $[0,c\lambda_0]$, $[-c\lambda_0,0]$ and $[-\lambda_0,-c\lambda_0]$, respectively. From the frequency constraint condition we can obtain the corresponding range of $\lambda_3$ (see Table \ref{table1}).
\begin{table}
\begin{center}

\begin{tabular}{|c|c|c|c|}
\hline
 ${\rm Case}$   &  $\lambda_2\in $ & $\lambda_1\in $ & $\lambda_3\in $    \\
\hline
$h_-h h$  & $ [-\lambda_0,-c\lambda_0]$  & $ [c\lambda_0, \lambda_0]$  &  $[\lambda_0, 2\lambda_0-c\lambda_0-l]$  \\
\hline
$h_-l h$  & $ [-\lambda_0,-c\lambda_0]$  & $ [0,c\lambda_0]$  &  $[\lambda_0, 2\lambda_0-l]$  \\
\hline
$h_-l_-h$ & $ [-\lambda_0,-c\lambda_0]$  & $ [-c\lambda_0,0]$ & $[\lambda_0+c\lambda_0-l, 2\lambda_0+c\lambda_0-l]$  \\
\hline
$h_-h_- h$ & $ [-\lambda_0,-c\lambda_0]$  & $ [-\lambda_0,-c\lambda_0]$ &$ [\lambda_0+2c\lambda_0-l, 3\lambda_0-l]$\\
\hline
\end{tabular}
\end{center}
\caption{$\lambda_2\in [-\lambda_0,-c\lambda_0]$}\label{table1}
\end{table}

{\it Case $h_-hh$.}  We decompose $\lambda_1, \lambda_2, \lambda_3$ by:
\begin{align}
\lambda_1\in [c\lambda_0,\lambda_0]= \bigcup_{j_1\geq 0}\lambda_0-I_{j_1},\ \ \lambda_2\in [-\lambda_0,-c\lambda_0]= \bigcup_{j_2\geq 0} -\lambda_0+I_{j_2},\nonumber\\
\lambda_3 \in [\lambda_0, (2-c)\lambda_0-l]= \bigcup_{j_3\geq 0} \lambda_0+I_{j_3}, \ \ j_1, j_2, j_3\leq \log_2 \langle\lambda_0\rangle.    \nonumber
\end{align}
From the frequency constraint condition \eqref{FCC}, we have
\begin{align}
2^{j_1}\approx 2^{j_2}+2^{j_3}.
\end{align}
It follows that $j_1\approx j_2\vee j_3$. When $j_1\approx j_2\geq j_3$, we can get the result by using the similar technique as that used in Case 1 of Step 1. When $j_1\approx j_3\geq j_2$, we just need to exchange the positions of $j_2$ and $j_3$ and use the similar way to obtain our conclusion. We omit the details.

{\it Case $h_-lh$.}  We decompose $\lambda_1, \lambda_2, \lambda_3$ by:
\begin{align}
\lambda_1\in [0,c\lambda_0]= \bigcup_{j_1\geq 0} I_{j_1},\ \ \lambda_2\in [-\lambda_0,-c\lambda_0]= \bigcup_{j_2\geq 0} -\lambda_0+I_{j_2},\nonumber\\
\lambda_3 \in [\lambda_0, 2\lambda_0-l]= \bigcup_{j_3\geq 0} \lambda_0+I_{j_3}, \ \ j_1, j_2, j_3\leq \log_2 \langle\lambda_0\rangle.    \nonumber
\end{align}
From the frequency constraint condition \eqref{FCC}, we have
\begin{align}
2^{j_1}+2^{j_2}+2^{j_3}\approx \lambda_0.
\end{align}
By DMCC \eqref{DMCC} the highest dispersion modulation satisfies
\begin{align}
\max_{0\leq k \leq 3} |\xi^3_k - \tau_k|  \gtrsim \langle\lambda_0\rangle^2\cdot 2^{j_3}.
\end{align}

If $v_{\lambda_0}$ has the highest dispersion modulation, we have dispersion modulation decay to $v_{\lambda_0}$. For $u_{I_{j_1}}$ and $u_{\lambda_0+I_{j_3}}$, we have $|\lambda_0+2^{j_3}+2^{j_1}|\gtrsim \langle\lambda_0\rangle$ and $|\lambda_0+2^{j_3}-2^{j_1}|\gtrsim \langle\lambda_0\rangle$. Thus we can use bilinear estimate \eqref{bilinear3} to $u_{I_{j_1}}u_{\lambda_0+I_{j_3}}$. To be specific, we have
\begin{align}\label{step2e}
 &\sum_{j_1, j_2, j_3\leq \log_2 \langle\lambda_0\rangle} \langle \lambda_0\rangle^{1/4}\int_{[0,T]\times\mathbb{R}} |\overline{v_{\lambda_0}} u_{I_{j_1}}u_{-\lambda_0+I_{j_2}}\partial_{x}u_{\lambda_0+I_{j_3}}| \, dxdt\nonumber\\
 \lesssim & \sum_{j_1, j_2, j_3\leq \log_2 \langle\lambda_0\rangle} \langle \lambda_0\rangle^{5/4}\|\overline{v_{\lambda_0}}\|_{L^2_tL^{\infty}_x} \| u_{-\lambda_0+I_{j_2}}\|_{L^\infty_tL^2_x} \|u_{I_{j_1}}u_{\lambda_0+I_{j_3}}\|_{L^2_{x,t}}\nonumber\\
\lesssim & \sum_{j_1, j_2, j_3\leq \log_2 \langle\lambda_0\rangle} \langle \lambda_0\rangle^{5/4}\langle \lambda_0\rangle^{-1}(2^{j_3})^{-1/2}\|v_{\lambda_0}\|_{V^2_A} \|u_{-\lambda_0+I_{j_2}}\|_{V^2_A}\nonumber\\
&\quad \times T^{\varepsilon/4}\langle \lambda_0\rangle^{-1+2\varepsilon}\|u_{I_{j_1}}\|_{V^2_A}\|u_{\lambda_0+I_{j_3}}\|_{V^2_A}\nonumber\\
\lesssim & T^{\varepsilon/4}\sum_{j_1, j_2, j_3\leq \log_2 \langle\lambda_0\rangle} \langle \lambda_0\rangle^{-5/4+2\varepsilon} (2^{j_1})^{1/4}(2^{j_2})^{1/2}\|v_{\lambda_0}\|_{V^2_A}\|u\|^3_{X^{1/4}_{\infty,A}}\nonumber\\
\lesssim &T^{\varepsilon/4}\langle \lambda_0\rangle^{-1/2+2\varepsilon}\log_2 \langle\lambda_0\rangle  \|v_{\lambda_0}\|_{V^2_A}\|u\|^3_{X^{1/4}_{\infty,A}}\nonumber\\
\lesssim &T^{\varepsilon/4} \|v^{(\lambda_0)}\|_{V^2_A}\|u\|^3_{X^{1/4}_{\infty,A}}.
\end{align}

If $u_{-\lambda_0+I_{j_2}}$ has the highest dispersion modulation,  we take $L^\infty_{x,t}$, $L^2_{x,t}$ and $L^2_{x,t}$ norms to $v_{\lambda_0}$, $u_{-\lambda_0+I_{j_2}}$, and $u_{I_{j_1}}u_{\lambda_0+I_{j_3}}$, respectively. Then applying the dispersion modulation decay estimate \eqref{dispersiondecay} to $u_{-\lambda_0+I_{j_2}}$  and the bilinear estimate \eqref{bilinear3} to $u_{I_{j_1}}u_{\lambda_0+I_{j_3}}$, we can get the desired conclusion.

If $u_{I_{j_1}}$ has the highest dispersion modulation, we have  dispersion modulation decay to $u_{I_{j_1}}$. For $u_{-\lambda_0+I_{j_2}}$ and $u_{\lambda_0+I_{j_3}}$, we have $|\lambda_0+2^{j_3}-\lambda_0+2^{j_2}|\gtrsim (2^{j_3}+2^{j_2})\approx \lambda_0-2^{j_1}\gtrsim \langle\lambda_0\rangle$ and $|\lambda_0+2^{j_3}+\lambda_0-2^{j_2}|\gtrsim \langle\lambda_0\rangle$. Thus we can use bilinear estimate \eqref{bilinear3} to $u_{-\lambda_0+I_{j_2}}u_{\lambda_0+I_{j_3}}$. Therefore, we have
\begin{align}\label{step2d}
&\sum_{j_1, j_2,j_3\leq \log_2 \langle\lambda_0\rangle} \langle \lambda_0\rangle^{1/4}\int_{[0,T]\times\mathbb{R}} |\overline{v_{\lambda_0}} u_{I_{j_1}}u_{-\lambda_0+I_{j_2}}\partial_{x}u_{\lambda_0+I_{j_3}}| \, dxdt\nonumber\\
&\lesssim \sum_{j_1, j_2, j_3\leq \log_2 \langle\lambda_0\rangle} \langle \lambda_0\rangle^{5/4}\|\overline{v_{\lambda_0}}\|_{L^\infty_{x,t}} \|u_{I_{j_1}}\|_{L^2_{x,t}}\|u_{-\lambda_0+I_{j_2}}u_{\lambda_0+I_{j_3}}\|_{L^2_{x,t}}\nonumber\\
&\lesssim \sum_{j_1, j_2, j_3\leq \log_2 \langle\lambda_0\rangle}\langle \lambda_0\rangle^{5/4}\|v_{\lambda_0}\|_{V^2_A} \langle \lambda_0\rangle^{-1}(2^{j_3})^{-1/2} \|u_{I_{j_1}}\|_{V^2_A}\nonumber\\
&\quad\times T^{\varepsilon/4}\langle \lambda_0\rangle^{-1+2\varepsilon}\|u_{-\lambda_0+I_{j_2}}\|_{V^2_A}\|u_{\lambda_0+I_{j_3}}\|_{V^2_A},
\end{align}
which is the same with the third line of \eqref{step2e}, so we omit the details.

If $u_{\lambda_0+I_{j_3}}$ has the highest dispersion modulation,  from the dispersion modulation decay \eqref{dispersiondecay} and $L^4$ estimate \eqref{lebesgue4a}, we have
\begin{align}\label{step2f}
 &\sum_{j_1, j_2, j_3\leq \log_2 \langle\lambda_0\rangle} \langle \lambda_0\rangle^{1/4}\int_{[0,T]\times\mathbb{R}} |\overline{v_{\lambda_0}} u_{I_{j_1}}u_{-\lambda_0+I_{j_2}}\partial_{x}u_{\lambda_0+I_{j_3}}| \, dxdt\nonumber\\
 \lesssim & \sum_{j_1, j_2, j_3\leq \log_2 \langle\lambda_0\rangle} \langle \lambda_0\rangle^{5/4}\|\overline{v_{\lambda_0}}\|_{L^\infty_{x,t}}\| u_{\lambda_0+I_{j_3}}\|_{L^2_{x,t}} \|u_{I_{j_1}}\|_{L^4_{x,t}}\|u_{-\lambda_0+I_{j_2}}\|_{L^4_{x,t}}\nonumber\\
\lesssim & \sum_{j_1, j_2, j_3\leq \log_2 \langle\lambda_0\rangle} \langle \lambda_0\rangle^{5/4}\|v_{\lambda_0}\|_{V^2_A} \langle \lambda_0\rangle^{-1}(2^{j_3})^{-1/2}(2^{j_3})^{1/2}\langle\lambda_0\rangle^{-1/4}\|u\|_{X^{1/4}_{\infty,A}}\nonumber\\
&\quad \times T^{\varepsilon/4} (2^{j_1})^{-1/8+\varepsilon}T^{\varepsilon/4} (2^{j_2})^{1/4+\varepsilon}\langle \lambda_0\rangle^{-3/8}\|u\|^2_{X^{1/4}_{\infty,A}}\nonumber\\
\lesssim & T^{\varepsilon/2} \sum_{j_1, j_2, j_3\leq \log_2 \langle\lambda_0\rangle} \langle \lambda_0\rangle^{-3/8}(2^{j_1})^{-1/8+\varepsilon}(2^{j_2})^{1/4+\varepsilon}\|v_{\lambda_0}\|_{V^2_A}\|u\|
^3_{X^{1/4}_{\infty,A}}.
\end{align}
Taking $0<\varepsilon<1/8$, the summation over $j_1$ is finite.  The summation over $j_2$ and $j_3$ can be controlled, so \eqref{step2f} is continued by
\begin{align}
\lesssim T^{\varepsilon/2}\langle \lambda_0\rangle^{-1/8+\varepsilon}\log_2\langle\lambda_0\rangle\|v_{\lambda_0}\|_{V^2_A}\|u\|^3_{X^{1/4}_{\infty,A}}
\lesssim  T^{\varepsilon/2}\|v^{(\lambda_0)}\|_{V^2_A}\|u\|^3_{X^{1/4}_{\infty,A}}.
\end{align}

{\it Case $h_-l_-h$.}  We decompose $\lambda_1, \lambda_2, \lambda_3$ in the following way:
\begin{gather}
\lambda_1\in [-c\lambda_0,0]= \bigcup_{j_1\geq 0}-I_{j_1},\ \ \lambda_2\in [-\lambda_0,-c\lambda_0]= \bigcup_{j_2\geq 0} -\lambda_0+I_{j_2},\ \ j_1,j_2\leq \log_2\langle\lambda_0\rangle, \nonumber\\
\lambda_3 \in [\lambda_0+c\lambda_0-l, 2\lambda_0+c\lambda_0-l]= \bigcup_{j_3\geq \log_2\langle\lambda_0\rangle-C} \lambda_0+I_{j_3}, \ \   j_3\leq \log_2\langle\lambda_0\rangle+1.  \nonumber
\end{gather}
From the frequency constraint condition \eqref{FCC}, we have
\begin{align}
2^{j_2}+2^{j_3}\approx \lambda_0+2^{j_1}.
\end{align}
It is easy to see that this case is similar to the above Case $h_-lh$, so the details are omitted.

{\it Case $h_-h_-h$.} We decompose $\lambda_1, \lambda_2, \lambda_3$ as follows:
\begin{gather}
\lambda_k\in[-\lambda_0,-c\lambda_0]=\bigcup_{j_k\geq 0} -\lambda_0+I_{j_k},\ \  j_k \leq \log_2\langle\lambda_0\rangle,\ \ k=1,2;  \\
\lambda_3\in[\lambda_0+2c\lambda_0-l, 3\lambda_0-l]=\bigcup_{j_3\geq \log_2(2c\lambda_0-l)} \lambda_0+I_{j_3},\ \ j_3\leq \log_2\langle\lambda_0\rangle+1.
\end{gather}
From the dispersion modulation constraint condition \eqref{DMCC}, we know the highest dispersion modulation satisfies
\begin{align}
\max_{0\leq k \leq 3} |\xi^3_k - \tau_k|  \gtrsim \langle\lambda_0\rangle^2\cdot 2^{j_3}.
\end{align}

If $v_{\lambda_0}$ has the highest dispersion modulation, we take dispersion modulation decay to $v_{\lambda_0}$. For $u_{-\lambda_0+I_{j_1}}$ and $u_{\lambda_0+I_{j_3}}$, we have $|\lambda_0+2^{j_3}-\lambda_0+2^{j_1}|\gtrsim 2^{j_1}$ and $|\lambda_0+2^{j_3}+\lambda_0-2^{j_1}|\gtrsim \langle\lambda_0\rangle$. Thus we can use bilinear estimate \eqref{bilinear3} to $u_{-\lambda_0+I_{j_1}}u_{\lambda_0+I_{j_3}}$. Thus we have
\begin{align}\label{step2g}
 &\sum_{j_1, j_2, j_3\leq \log_2 \langle\lambda_0\rangle+1} \langle \lambda_0\rangle^{1/4}\int_{[0,T]\times\mathbb{R}} |\overline{v_{\lambda_0}} u_{-\lambda_0+I_{j_1}}u_{-\lambda_0+I_{j_2}}\partial_{x}u_{\lambda_0+I_{j_3}}| \, dxdt\nonumber\\
 \lesssim & \sum_{j_1, j_2, j_3\leq \log_2 \langle\lambda_0\rangle+1} \langle \lambda_0\rangle^{5/4}\|\overline{v_{\lambda_0}}\|_{L^2_tL^{\infty}_x} \| u_{-\lambda_0+I_{j_2}}\|_{L^\infty_tL^2_x} \|u_{-\lambda_0+I_{j_1}}u_{\lambda_0+I_{j_3}}\|_{L^2_{x,t}}\nonumber\\
\lesssim & \sum_{j_1, j_2, j_3\leq \log_2 \langle\lambda_0\rangle+1} \langle \lambda_0\rangle^{5/4}\langle \lambda_0\rangle^{-1}(2^{j_3})^{-1/2}\|v_{\lambda_0}\|_{V^2_A} \|u_{-\lambda_0+I_{j_2}}\|_{V^2_A}\nonumber\\
&\quad \times T^{\varepsilon/4}\langle \lambda_0\rangle^{-1/2+\varepsilon}(2^{j_1})^{-1/2+\varepsilon}\|u_{-\lambda_0+I_{j_1}}\|_{V^2_A}\|u_{\lambda_0+I_{j_3}}\|_{V^2_A}\nonumber\\
\lesssim & T^{\varepsilon/4}\sum_{j_1, j_2, j_3\leq \log_2 \langle\lambda_0\rangle+1} \langle \lambda_0\rangle^{-1+\varepsilon} (2^{j_1})^{\varepsilon}(2^{j_2})^{1/2}\|v_{\lambda_0}\|_{V^2_A}\|u\|^3_{X^{1/4}_{\infty,A}}\nonumber\\
\lesssim &T^{\varepsilon/4}\langle \lambda_0\rangle^{-1/2+2\varepsilon}\log_2 \langle\lambda_0\rangle  \|v_{\lambda_0}\|_{V^2_A}\|u\|^3_{X^{1/4}_{\infty,A}}\nonumber\\
\lesssim &T^{\varepsilon/4} \|v^{(\lambda_0)}\|_{V^2_A}\|u\|^3_{X^{1/4}_{\infty,A}}.
\end{align}

If $u_{-\lambda_0+I_{j_2}}$ has the highest dispersion modulation, we take $L^\infty_{x,t}$, $L^2_{x,t}$ and $L^2_{x,t}$ norms to $v_{\lambda_0}$, $u_{-\lambda_0+I_{j_2}}$, and $u_{-\lambda_0+I_{j_1}}u_{\lambda_0+I_{j_3}}$, respectively. Then we can get the desired estimate by using the analogue technique. If $u_{-\lambda_0+I_{j_1}}$ has the highest dispersion modulation, due to the symmetry between $u_{-\lambda_0+I_{j_1}}$ and $u_{-\lambda_0+I_{j_2}}$, the estimate is similar so we omit the details.

If $u_{\lambda_0+I_{j_3}}$ has the highest dispersion modulation,  from the dispersion modulation decay \eqref{dispersiondecay} and $L^4$ estimate \eqref{lebesgue4a}, we have
\begin{align}\label{step2h}
 &\sum_{j_1, j_2, j_3\leq \log_2 \langle\lambda_0\rangle+1} \langle \lambda_0\rangle^{1/4}\int_{[0,T]\times\mathbb{R}} |\overline{v_{\lambda_0}} u_{-\lambda_0+I_{j_1}}u_{-\lambda_0+I_{j_2}}\partial_{x}u_{\lambda_0+I_{j_3}}| \, dxdt\nonumber\\
 \lesssim & \sum_{j_1, j_2, j_3\leq \log_2 \langle\lambda_0\rangle+1} \langle \lambda_0\rangle^{5/4}\|\overline{v_{\lambda_0}}\|_{L^\infty_{x,t}}\| u_{\lambda_0+I_{j_3}}\|_{L^2_{x,t}} \|u_{-\lambda_0+I_{j_1}}\|_{L^4_{x,t}}\|u_{-\lambda_0+I_{j_2}}\|_{L^4_{x,t}}\nonumber\\
\lesssim & \sum_{j_1, j_2, j_3\leq \log_2 \langle\lambda_0\rangle+1} \langle \lambda_0\rangle^{5/4}\|v_{\lambda_0}\|_{V^2_A} \langle \lambda_0\rangle^{-1}(2^{j_3})^{-1/2}(2^{j_3})^{1/2}\langle\lambda_0\rangle^{-1/4}\|u\|_{X^{1/4}
_{\infty,A}}\nonumber\\
&\quad \times T^{\varepsilon/4} (2^{j_1})^{1/4+\varepsilon}\langle \lambda_0\rangle^{-3/8}T^{\varepsilon/4} (2^{j_2})^{1/4+\varepsilon}\langle \lambda_0\rangle^{-3/8}\|u\|^2_{X^{1/4}_{\infty,A}}\nonumber\\
\lesssim & T^{\varepsilon/2} \sum_{j_1, j_2, j_3\leq \log_2 \langle\lambda_0\rangle+1} \langle \lambda_0\rangle^{-3/4}(2^{j_1})^{1/4+\varepsilon}(2^{j_2})^{1/4+\varepsilon}\|v_{\lambda_0}\|_{V^2_A}
\|u\|^3_{X^{1/4}_{\infty,A}}\nonumber\\
\lesssim & T^{\varepsilon/2}\langle \lambda_0\rangle^{-1/4+2\varepsilon}\log_2 \langle\lambda_0\rangle\|v_{\lambda_0}\|_{V^2_A}\|u\|^3_{X^{1/4}_{\infty,A}}\nonumber\\
\lesssim & T^{\varepsilon/2}\|v^{(\lambda_0)}\|_{V^2_A}\|u\|^3_{X^{1/4}_{\infty,A}}.\nonumber
\end{align}

{\bf Case 3.3:} $\lambda_2\in (-\infty,-\lambda_0]$. We consider $\lambda_1\in [c\lambda_0, \lambda_0]$, $[0,c\lambda_0]$, $[-\lambda_0,0]$ and $(-\infty,-\lambda_0]$, respectively. From the frequency constraint condition we can obtain the corresponding range of $\lambda_3$ (see Table \ref{table2}).
\begin{table}
\begin{center}

\begin{tabular}{|c|c|c|c|}
\hline
 ${\rm Case}$   &  $\lambda_2\in $ & $\lambda_1\in $ & $\lambda_3\in $    \\
\hline
$2h_-hh$  & $ (-\infty,-\lambda_0]$  & $ [c\lambda_0,\lambda_0]$  &  $[\lambda_0, +\infty]$  \\
\hline
$2h_-lh$  & $ (-\infty,-\lambda_0]$  & $ [0,c\lambda_0]$  &  $[2\lambda_0-c\lambda_0-l, 2\lambda_0]$  \\
\hline
$2h_-lh2$  & $ (-\infty,-\lambda_0]$  & $ [0,c\lambda_0]$  &  $[2\lambda_0,\infty]$  \\
\hline
$2h_-l_-h$ & $ (-\infty,-\lambda_0]$  & $ [-\lambda_0,0]$ &$ [2\lambda_0-l, +\infty]$\\
\hline
$2h_-h_-h$ & $ (-\infty,-\lambda_0]$  & $ (-\infty,-\lambda_0]$ &$ [3\lambda_0-l, +\infty]$\\
\hline
\end{tabular}
\end{center}
\caption{$\lambda_2\in (-\infty,-\lambda_0]$}\label{table2}
\end{table}

 Case $2h_-hh$.  We decompose $\lambda_1, \lambda_2, \lambda_3$ in the following way:
\begin{gather*}
\lambda_1\in[c\lambda_0,\lambda_0]=\bigcup_{j_1\geq 0} \lambda_0-I_{j_1},\ \  j_1 \leq \log_2\langle\lambda_0\rangle;  \\
\lambda_2\in(-\infty, -\lambda_0]=\bigcup_{j_2\geq 0} -\lambda_0-I_{j_2},\ \ \lambda_3\in[\lambda_0,+\infty)=\bigcup_{j_3\geq 0} \lambda_0+I_{j_3}.
\end{gather*}
From the frequency constraint condition \eqref{FCC}, we know that
\begin{align}
2^{j_3}\approx 2^{j_1}+2^{j_2} \Rightarrow j_3\approx j_1\vee j_2.\nonumber
\end{align}
From the dispersion modulation constraint condition \eqref{DMCC}, we have the highest dispersion modulation satisfies
\begin{align}
\max_{0\leq k \leq 3} |\xi^3_k - \tau_k|  \gtrsim (\langle\lambda_0\rangle+2^{j_2})\cdot 2^{j_1}\cdot 2^{j_3}.
\end{align}

If $v_{\lambda_0}$ has the highest dispersion modulation, from the dispersion modulation decay \eqref{dispersiondecay} and $L^4$ estimate \eqref{lebesgue4a}, we have
\begin{align*}
 &\sum_{j_1,j_2,j_3\geq 0} \langle \lambda_0\rangle^{1/4}\int_{[0,T]\times\mathbb{R}} |\overline{v_{\lambda_0}} u_{\lambda_0-I_{j_1}}u_{-\lambda_0-I_{j_2}}\partial_{x}u_{\lambda_0+I_{j_3}}| \, dxdt\nonumber\\
 \lesssim & \sum_{j_1,j_2,j_3\geq 0} \langle \lambda_0\rangle^{1/4}(\langle \lambda_0\rangle+2^{j_3})\|\overline{v_{\lambda_0}}\|_{L^2_tL^\infty_x}\| u_{\lambda_0-I_{j_1}}\|_{L^\infty_tL^2_x} \|u_{-\lambda_0-I_{j_2}}\|_{L^4_{x,t}}\|u_{\lambda_0+I_{j_3}}\|_{L^4_{x,t}}\nonumber\\
\lesssim & \sum_{j_1,j_2,j_3\geq 0} (\langle \lambda_0\rangle+2^{j_3})(\langle \lambda_0\rangle+2^{j_2})^{-1/2}(2^{j_3})^{-1/2}\|v_{\lambda_0}\|_{V^2_A} \|u\|_{X^{1/4}_{\infty,A}}\nonumber\\
&\quad \times T^{\varepsilon/2} (2^{j_2})^{1/4+\varepsilon}(\langle \lambda_0\rangle+2^{j_2})^{-3/8} (2^{j_3})^{1/4+\varepsilon}(\langle \lambda_0\rangle+2^{j_3})^{-3/8}\|u\|^2_{X^{1/4}_{\infty,A}}\nonumber\\
\lesssim & T^{\varepsilon/2} \sum_{j_1, j_2, j_3\geq 0} (\langle \lambda_0\rangle+2^{j_3})^{5/8}(\langle \lambda_0\rangle+2^{j_2})^{-7/8}(2^{j_3})^{-1/4+\varepsilon}(2^{j_2})^{1/4+\varepsilon}\|v_{\lambda_0}\|_{V^2_A}
\|u\|^3_{X^{1/4}_{\infty,A}}.
\end{align*}
When $0\leq j_1\leq j_2\approx j_3$, the summation in above inequality becomes
\begin{align}
\sum_{j_3\geq 0} (\langle \lambda_0\rangle+2^{j_3})^{-1/4}(2^{j_3})^{2\varepsilon}\cdot j_3 \lesssim 1.
\end{align}
When $0\leq j_2\leq j_1\approx j_3$, noticing that $j_1\leq \log_2\langle \lambda_0\rangle$, we can know that the summation  satisfies
\begin{align}
\sum_{0\leq j_2\leq j_3\leq \log_2\langle \lambda_0\rangle} \langle \lambda_0\rangle^{-1/4}(2^{j_3})^{-1/4+\varepsilon}(2^{j_2})^{1/4+\varepsilon}\cdot j_3 \lesssim 1.
\end{align}

If $u_{\lambda_0-I_{j_1}}$ has the highest dispersion modulation, we take $L^\infty_{x,t}$, $L^2_{x,t}$, $L^4_{x,t}$ and $L^4_{x,t}$ norms to $v_{\lambda_0}$, $u_{\lambda_0-I_{j_1}}$, $u_{-\lambda_0-I_{j_2}}$ and $u_{\lambda_0+I_{j_3}}$, respectively. Then we can get the desired estimate by a similar way.

If $u_{-\lambda_0-I_{j_2}}$ has the highest dispersion modulation, we take the dispersion modulation decay estimate to $u_{-\lambda_0-I_{j_2}}$. For $u_{\lambda_0-I_{j_1}}$ and $u_{\lambda_0+I_{j_3}}$, we have $|\lambda_0+2^{j_3}-\lambda_0+2^{j_1}|\gtrsim 2^{j_3}$ and $|\lambda_0+2^{j_3}+\lambda_0-2^{j_1}|\gtrsim ( \langle\lambda_0\rangle+2^{j_3})$. Thus we can use the bilinear estimate \eqref{bilinear3} to $u_{\lambda_0-I_{j_1}}u_{\lambda_0+I_{j_3}}$. Thus we have
\begin{align*}
 &\sum_{j_1, j_2, j_3\geq 0} \langle \lambda_0\rangle^{1/4}\int_{[0,T]\times\mathbb{R}} |\overline{v_{\lambda_0}} u_{\lambda_0-I_{j_1}}u_{-\lambda_0-I_{j_2}}\partial_{x}u_{\lambda_0+I_{j_3}}| \, dxdt\nonumber\\
 \lesssim & \sum_{j_1, j_2, j_3\geq 0} \langle \lambda_0\rangle^{1/4}\|\overline{v_{\lambda_0}}\|_{L^\infty_tL^{\infty}_x}\| u_{-\lambda_0-I_{j_2}}\|_{L^2_tL^2_x} \|u_{\lambda_0-I_{j_1}}\partial_x u_{\lambda_0+I_{j_3}}\|_{L^2_{x,t}}\nonumber\\
\lesssim & \sum_{j_1, j_2, j_3\geq 0} \langle \lambda_0\rangle^{1/4}(\langle \lambda_0\rangle+2^{j_3})(\langle \lambda_0\rangle+2^{j_2})^{-1/2}(2^{j_1})^{-1/2}(2^{j_3})^{-1/2}\|v_{\lambda_0}\|_{V^2_A} \|u_{-\lambda_0-I_{j_2}}\|_{V^2_A}\nonumber\\
&\quad \times T^{\varepsilon/4}(\langle \lambda_0\rangle+2^{j_3})^{-1/2+\varepsilon}(2^{j_3})^{-1/2+\varepsilon}\|u_{\lambda_0-I_{j_1}}\|_{V^2_A}
\|u_{\lambda_0+I_{j_3}}\|_{V^2_A}\nonumber\\
\lesssim & T^{\varepsilon/4}\sum_{j_1, j_2, j_3\geq 0} (\langle \lambda_0\rangle+2^{j_3})^{1/4+\varepsilon}(\langle \lambda_0\rangle+2^{j_2})^{-3/4}(2^{j_3})^{-1/2+\varepsilon}(2^{j_2})^{1/2}\|v_{\lambda_0}\|_{V^2_A}
\|u\|^3_{X^{1/4}_{\infty,A}}.
\end{align*}
If $0\leq j_1\leq j_2\approx j_3$, the summation in above inequality becomes
\begin{align*}
\sum_{j_3\geq 0} (\langle \lambda_0\rangle+2^{j_3})^{-1/2+\varepsilon}(2^{j_3})^{\varepsilon}\cdot j_3 \lesssim 1.
\end{align*}
If $0\leq j_2\leq j_1\approx j_3$, recalling that $j_1\leq \log_2\langle \lambda_0\rangle$, we can get the summation  satisfying
\begin{align*}
\sum_{0\leq j_2\leq j_3\leq \log_2\langle \lambda_0\rangle} \langle \lambda_0\rangle^{-1/2+\varepsilon}(2^{j_3})^{-1/2+\varepsilon}(2^{j_2})^{1/2}\cdot j_3 \lesssim 1.
\end{align*}

If $u_{\lambda_0+I_{j_3}}$ attains the highest dispersion modulation, noticing that for $u_{\lambda_0-I_{j_1}}$ and $u_{-\lambda_0-I_{j_2}}$, we have $|\lambda_0+2^{j_2}-\lambda_0+2^{j_1}|\gtrsim 2^{j_2}$ and $|\lambda_0+2^{j_2}+\lambda_0-2^{j_1}|\gtrsim ( \langle\lambda_0\rangle+2^{j_2})$. We can use the bilinear estimate \eqref{bilinear3} to $u_{\lambda_0-I_{j_1}}u_{-\lambda_0-I_{j_3}}$ to get our result by using the same way as above.

 Case $2h_-lh$. From (FCC) \eqref{FCC}, we see that $\lambda_2\in [-\lambda_0-c\lambda_0-l,-\lambda_0]$.  We decompose $\lambda_1, \lambda_2, \lambda_3$ in a dyadic way:
\begin{align}
\lambda_1\in [0,c\lambda_0]= \bigcup_{j_1\geq 0} I_{j_1},\ \ \lambda_2\in [-\lambda_0-c\lambda_0-l,-\lambda_0]= \bigcup_{j_2\geq 0} -\lambda_0-I_{j_2},\nonumber\\
\lambda_3 \in [2\lambda_0-c\lambda_0-l, 2\lambda_0]= \bigcup_{j_3\geq 0} \lambda_0+I_{j_3}, \ \ j_1, j_2, j_3\leq \log_2 \langle\lambda_0\rangle.    \nonumber
\end{align}
From the frequency constraint condition \eqref{FCC}, we have
\begin{align}
2^{j_1}+2^{j_3}-2^{j_2}\approx \lambda_0.\nonumber
\end{align}
By DMCC \eqref{DMCC} the highest dispersion modulation satisfies
\begin{align}
\max_{0\leq k \leq 3} |\xi^3_k - \tau_k|  \gtrsim \langle\lambda_0\rangle^2\cdot 2^{j_3}.\nonumber
\end{align}
Therefore, the approach to this case is similar to Case $h_-lh$, and we omit it.

Case $2h_-lh2$. We decompose $\lambda_1, \lambda_2, \lambda_3$ in the following way:
\begin{gather}
\lambda_1\in [0,c\lambda_0]= \bigcup_{j_1\geq 0} I_{j_1},\ \ j_1\leq\log_2 \langle\lambda_0\rangle; \nonumber\\
\lambda_2\in [-\infty,-\lambda_0]= \bigcup_{j_2\geq 0} -\lambda_0-I_{j_2},\ \ \lambda_3 \in [2\lambda_0,+\infty]= \bigcup_{j_3\geq 0} 2\lambda_0+I_{j_3}. \nonumber
\end{gather}
From the frequency constraint condition \eqref{FCC}, we have
\begin{align}
2^{j_2}\approx2^{j_1}+2^{j_3},\ \ \ {\rm i.e.}\ \ \  j_2\approx j_1\vee j_3.
\end{align}
By DMCC \eqref{DMCC} the highest dispersion modulation satisfies
\begin{align}
\max_{0\leq k \leq 3} |\xi^3_k - \tau_k|  \gtrsim \langle\lambda_0\rangle\cdot(\langle\lambda_0\rangle+2^{j_2})\cdot( \langle\lambda_0\rangle+2^{j_3}).
\end{align}

If $v_{\lambda_0}$ has the highest dispersion modulation, we take dispersion modulation decay to $v_{\lambda_0}$. For $u_{I_{j_1}}$ and $u_{2\lambda_0+I_{j_3}}$, we have $|2\lambda_0+2^{j_3}\pm2^{j_1}|\gtrsim (\langle\lambda_0\rangle+2^{j_3})$. Thus we can use bilinear estimate \eqref{bilinear3} to $u_{I_{j_1}}u_{2\lambda_0+I_{j_3}}$. Specifically, we have
\begin{align}\label{step2i}
 &\sum_{j_2\approx j_1\vee j_3} \langle \lambda_0\rangle^{1/4}\int_{[0,T]\times\mathbb{R}} |\overline{v_{\lambda_0}} u_{I_{j_1}}u_{-\lambda_0-I_{j_2}}\partial_{x}u_{2\lambda_0+I_{j_3}}| \, dxdt\nonumber\\
 \lesssim & \sum_{j_2\approx j_1\vee j_3} \langle \lambda_0\rangle^{1/4}(\langle\lambda_0\rangle+2^{j_3})\|\overline{v_{\lambda_0}}\|_{L^2_tL^{\infty}_x} \| u_{-\lambda_0-I_{j_2}}\|_{L^\infty_tL^2_x} \|u_{I_{j_1}}u_{2\lambda_0+I_{j_3}}\|_{L^2_{x,t}}\nonumber\\
\lesssim & \sum_{j_2\approx j_1\vee j_3} \langle \lambda_0\rangle^{1/4}\langle \lambda_0\rangle^{-1/2}(\langle\lambda_0\rangle+2^{j_2})^{-1/2}( \langle\lambda_0\rangle+2^{j_3})^{1/2}\|v_{\lambda_0}\|_{V^2_A} \|u_{-\lambda_0-I_{j_2}}\|_{V^2_A}\nonumber\\
&\quad \times T^{\varepsilon/4}(\langle\lambda_0\rangle+2^{j_3})^{-1+2\varepsilon}\|u_{I_{j_1}}\|_{V^2_A}\|u_{2\lambda_0+I_{j_3}}\|_{V^2_A}\nonumber\\
\lesssim & T^{\varepsilon/4}\sum_{j_2\approx j_1\vee j_3} \langle \lambda_0\rangle^{-1/4} ( \langle\lambda_0\rangle+2^{j_3})^{-3/4+2\varepsilon} ( \langle\lambda_0\rangle+2^{j_2})^{-3/4}\nonumber\\
&\quad\quad\quad \times (2^{j_1})^{1/4}(2^{j_2})^{1/2}(2^{j_3})^{1/2}\|v_{\lambda_0}\|
_{V^2_A}\|u\|^3_{X^{1/4}_{\infty,A}}\nonumber\\
\lesssim &T^{\varepsilon/4} \|v^{(\lambda_0)}\|_{V^2_A}\|u\|^3_{X^{1/4}_{\infty,A}},
\end{align}
where the last inequality is by summing over $j_1$, $j_2$ and $j_3$. Indeed we have the following estimates:
\begin{gather}
 \sum_{ j_1\leq \log_2 \langle\lambda_0\rangle }(2^{j_1})^{1/4}\lesssim \langle\lambda_0\rangle^{1/4}; \ \
 \sum_{j_2\geq 0}( \langle\lambda_0\rangle+2^{j_2})^{-3/4}(2^{j_2})^{1/2}\leq  \sum_{j_2\geq 0}(2^{j_2})^{-1/4}\lesssim 1;\nonumber\\
 \sum_{j_3\geq 0}( \langle\lambda_0\rangle+2^{j_3})^{-3/4+2\varepsilon}(2^{j_3})^{1/2}\leq  \sum_{j_3\geq 0}(2^{j_3})^{-1/4+2\varepsilon}\lesssim 1, \ \ 0<\varepsilon<1/8.\nonumber
\end{gather}

If $u_{-\lambda_0-I_{j_2}}$ has the highest dispersion modulation, we take $L^\infty_{x,t}$, $L^2_{x,t}$ and $L^2_{x,t}$ norms to $v_{\lambda_0}$, $u_{-\lambda_0-I_{j_2}}$, and $u_{I_{j_1}}u_{2\lambda_0+I_{j_3}}$, respectively. Then it will be same with \eqref{step2i}.

If $u_{2\lambda_0+I_{j_3}}$ has the highest dispersion modulation, we take dispersion modulation decay to $u_{2\lambda_0+I_{j_3}}$. For $u_{I_{j_1}}$ and $u_{-\lambda_0-I_{j_2}}$, we have $|\lambda_0+2^{j_2}+2^{j_1}|\gtrsim (\langle\lambda_0\rangle+2^{j_2})$ and $|\lambda_0+2^{j_2}-2^{j_1}|\gtrsim 2^{j_2}$. Thus we can use bilinear estimate \eqref{bilinear3} to $u_{I_{j_1}}u_{-\lambda_0-I_{j_3}}$. To be specific, we have
\begin{align*}
 &\sum_{j_2\approx j_1\vee j_3} \langle \lambda_0\rangle^{1/4}\int_{[0,T]\times\mathbb{R}} |\overline{v_{\lambda_0}} u_{I_{j_1}}u_{-\lambda_0-I_{j_2}}\partial_{x}u_{2\lambda_0+I_{j_3}}| \, dxdt\nonumber\\
 \lesssim & \sum_{j_2\approx j_1\vee j_3} \langle \lambda_0\rangle^{1/4}(\langle\lambda_0\rangle+2^{j_3})\|\overline{v_{\lambda_0}}\|_{L^\infty_{x,t}} \| u_{2\lambda_0+I_{j_3}}\|_{L^2_{x,t}} \|u_{I_{j_1}}u_{-\lambda_0-I_{j_2}}\|_{L^2_{x,t}}\nonumber\\
\lesssim & \sum_{j_2\approx j_1\vee j_3} \langle \lambda_0\rangle^{1/4}\|v_{\lambda_0}\|_{V^2_A} \langle \lambda_0\rangle^{-1/2}(\langle\lambda_0\rangle+2^{j_2})^{-1/2}( \langle\lambda_0\rangle+2^{j_3})^{1/2}\|u_{2\lambda_0+I_{j_3}}\|_{V^2_A}\nonumber\\
&\quad \times T^{\varepsilon/4}(\langle\lambda_0\rangle+2^{j_2})^{-1/2+\varepsilon}(2^{j_2})^{-1/2+\varepsilon}\|u_{I_{j_1}}
\|_{V^2_A}\|u_{-\lambda_0-I_{j_2}}\|_{V^2_A}\nonumber\\
\lesssim & T^{\varepsilon/4}\sum_{j_2\approx j_1\vee j_3} \langle \lambda_0\rangle^{-1/4} ( \langle\lambda_0\rangle+2^{j_3})^{1/4} ( \langle\lambda_0\rangle+2^{j_2})^{-5/4+\varepsilon}\nonumber\\
&\quad\quad\times(2^{j_1})^{1/4}(2^{j_2})^{\varepsilon}(2^{j_3})^{1/2}
\|v_{\lambda_0}\|_{V^2_A}\|u\|^3_{X^{1/4}_{\infty,A}}\nonumber\\
\lesssim &T^{\varepsilon/4} \|v^{(\lambda_0)}\|_{V^2_A}\|u\|^3_{X^{1/4}_{\infty,A}},
\end{align*}
where the last inequality is by summing over $j_1$, $j_2$ and $j_3$ in order.

If $u_{I_{j_1}}$ has the highest dispersion modulation,  from the dispersion modulation decay \eqref{dispersiondecay} and $L^4$ estimate \eqref{lebesgue4a}, we have
\begin{align}
 &\sum_{j_2\approx j_1\vee j_3} \langle \lambda_0\rangle^{1/4}\int_{[0,T]\times\mathbb{R}} |\overline{v_{\lambda_0}} u_{I_{j_1}}u_{-\lambda_0-I_{j_2}}\partial_{x}u_{2\lambda_0+I_{j_3}}| \, dxdt\nonumber\\
 \lesssim & \sum_{j_2\approx j_1\vee j_3} \langle \lambda_0\rangle^{1/4}(\langle\lambda_0\rangle+2^{j_3})\|\overline{v_{\lambda_0}}\|_{L^\infty_{x,t}}\| u_{I_{j_1}}\|_{L^2_{x,t}} \|u_{-\lambda_0-I_{j_2}}\|_{L^4_{x,t}}\|u_{2\lambda_0+I_{j_3}}\|_{L^4_{x,t}}\nonumber\\
\lesssim & \sum_{j_2\approx j_1\vee j_3} \langle \lambda_0\rangle^{1/4}\|v_{\lambda_0}\|_{V^2_A}  \langle \lambda_0\rangle^{-1/2}(\langle\lambda_0\rangle+2^{j_2})^{-1/2}( \langle\lambda_0\rangle+2^{j_3})^{1/2}(2^{j_1})^{1/4}\|u\|_{X^{1/4}_{\infty,A}}\nonumber\\
&\quad \times T^{\varepsilon/4} (2^{j_2})^{1/4+\varepsilon}(\langle\lambda_0\rangle+2^{j_2})^{-3/8}T^{\varepsilon/4} (2^{j_3})^{1/4+\varepsilon}(\langle\lambda_0\rangle+2^{j_3})^{-3/8}\|u\|^2_{X^{1/4}_{\infty,A}}\nonumber\\
\lesssim & T^{\varepsilon/2} \sum_{j_2\approx j_1\vee j_3} \langle \lambda_0\rangle^{-1/4}(2^{j_1})^{1/4}( \langle\lambda_0\rangle+2^{j_3})^{1/8} ( \langle\lambda_0\rangle+2^{j_2})^{-7/8}\nonumber\\
&\quad\quad\times(2^{j_2})^{1/4+\varepsilon}(2^{j_3})^{1/4+\varepsilon}\|v_{\lambda_0}
\|_{V^2_A}\|u\|^3_{X^{1/4}_{\infty,A}}\nonumber\\
\lesssim & T^{\varepsilon/2}\|v^{(\lambda_0)}\|_{V^2_A}\|u\|^3_{X^{1/4}_{\infty,A}},\nonumber
\end{align}
where the last inequality is by summing over $j_1$, $j_2$ and $j_3$ in order and noticing the condition $j_1\leq \log_2\langle\lambda_0\rangle$, $j_3\leq j_2$.

Case $2h_-l_-h$. We decompose $\lambda_1, \lambda_2, \lambda_3$ as follows:
\begin{gather}
\lambda_1\in [-\lambda_0,0]= \bigcup_{j_1\geq 0}-I_{j_1},\ \ j_1\leq\log_2 \langle\lambda_0\rangle; \nonumber\\
\lambda_2\in [-\infty,-\lambda_0]= \bigcup_{j_2\geq 0} -\lambda_0-I_{j_2},\ \ \lambda_3 \in [2\lambda_0-l,+\infty]= \bigcup_{j_3\geq -1} 2\lambda_0+I_{j_3}. \nonumber
\end{gather}
From the frequency constraint condition \eqref{FCC}, we have
\begin{align}
2^{j_3}\approx2^{j_1}+2^{j_2},\ \ \ {\rm i.e.}\ \ \  j_3\approx j_1\vee j_2.
\end{align}
If $j_3\approx j_2\geq j_1$, the method of this case will be same with Case $2h_-lh2$. If $j_3\approx j_1\geq j_2$, it is to say that $0\leq j_2 \leq j_3\approx j_1\leq \log_2 \langle\lambda_0\rangle$ holds, which can also ensure the convergence of the summation in Case $2h_-lh2$. Therefore, the details are omitted.

Case $2h_-h_-h$. We decompose $\lambda_1, \lambda_2, \lambda_3$ in the following way:
\begin{gather}
\lambda_k\in [-\infty, -\lambda_0]= \bigcup_{j_k\geq 0} -\lambda_0-I_{j_k}, k=1,2;\ \ \lambda_3 \in [3\lambda_0-l,+\infty]= \bigcup_{j_3\geq -1} 3\lambda_0+I_{j_3}. \nonumber
\end{gather}
From the frequency constraint condition \eqref{FCC} and $\lambda_1\geq\lambda_2$, we have
\begin{align}
2^{j_3}\approx2^{j_1}+2^{j_2},\  j_1\leq j_2\ \ \ {\rm i.e.}\ \ \  j_3\approx j_2\geq j_1.
\end{align}
By DMCC \eqref{DMCC} the highest dispersion modulation satisfies
\begin{align}
\max_{0\leq k \leq 3} |\xi^3_k - \tau_k|  \gtrsim (\langle\lambda_0\rangle+2^{j_1})\cdot(\langle\lambda_0\rangle+2^{j_2})\cdot( \langle\lambda_0\rangle+2^{j_3}).
\end{align}

If $v_{\lambda_0}$ has the highest dispersion modulation, from the dispersion modulation decay \eqref{dispersiondecay}, $L^4$ estimate \eqref{lebesgue4a}, and Lemma \ref{V2toX}, we have
\begin{align}\label{step2j}
 &\sum_{j_3\approx j_2\geq j_1} \langle \lambda_0\rangle^{1/4}\int_{[0,T]\times\mathbb{R}} |\overline{v_{\lambda_0}} u_{-\lambda_0-I_{j_1}}u_{-\lambda_0-I_{j_2}}\partial_{x}u_{3\lambda_0+I_{j_3}}| \, dxdt\nonumber\\
 \lesssim & \sum_{j_3\approx j_2\geq j_1} \langle \lambda_0\rangle^{1/4}(\langle\lambda_0\rangle+2^{j_3})\|\overline{v_{\lambda_0}}\|_{L^2_tL^{\infty}_x} \| u_{-\lambda_0-I_{j_1}}\|_{L^\infty_tL^2_x} \|u_{-\lambda_0-I_{j_2}}\|_{L^4_{x,t}}\|u_{3\lambda_0+I_{j_3}}\|_{L^4_{x,t}}\nonumber\\
\lesssim & \sum_{j_3\approx j_2\geq j_1} \langle \lambda_0\rangle^{1/4}(\langle\lambda_0\rangle+2^{j_1})^{-1/2}
(\langle\lambda_0\rangle+2^{j_2})^{-1/2}( \langle\lambda_0\rangle+2^{j_3})^{1/2}\|v_{\lambda_0}\|_{V^2_A}
(2^{j_1})^{1/2}\nonumber\\
&\quad \times  (\langle\lambda_0\rangle+2^{j_1})^{-1/4} T^{\varepsilon/2}(2^{j_2})^{1/4+\varepsilon}(\langle\lambda_0\rangle+2^{j_2})^{-3/8}
(2^{j_3})^{1/4+\varepsilon} (\langle\lambda_0\rangle+2^{j_3})^{-3/8}\|u\|^3_{X^{1/4}_{\infty,A}}\nonumber\\
\lesssim & T^{\varepsilon/2}\sum_{j_3\geq j_1} \langle \lambda_0\rangle^{1/4} ( \langle\lambda_0\rangle+2^{j_3})^{-3/4}(2^{j_3})^{1/2+2\varepsilon} ( \langle\lambda_0\rangle+2^{j_1})^{-3/4}(2^{j_1})^{1/2}\|v_{\lambda_0}\|_{V^2_A}\|u\|^3_{X^{1/4}_{\infty,A}}\nonumber\\
\lesssim & T^{\varepsilon/2}\bigg(\sum_{j_3\geq 0} (2^{j_3})^{-1/4+2\varepsilon}\cdot j_3\bigg) \|v_{\lambda_0}\|_{V^2_A}\|u\|^3_{X^{1/4}_{\infty,A}}\nonumber\\
\lesssim &T^{\varepsilon/2} \|v^{(\lambda_0)}\|_{V^2_A}\|u\|^3_{X^{1/4}_{\infty,A}}.
\end{align}

If $u_{-\lambda_0-I_{j_1}}$ has the highest dispersion modulation, we take $L^\infty_{x,t}$, $L^2_{x,t}$, $L^4_{x,t}$ and $L^4_{x,t}$ norms to $v_{\lambda_0}$, $u_{-\lambda_0-I_{j_1}}$, $u_{-\lambda_0-I_{j_2}}$ and $u_{3\lambda_0+I_{j_3}}$, respectively. Then it will be same with \eqref{step2j}.

If $u_{3\lambda_0+I_{j_3}}$ has the highest dispersion modulation, we divide the left hand side of \eqref{trilinear4} into two terms.
\begin{align}\label{step2m}
 &\sum_{j_3\approx j_2\geq j_1} \langle \lambda_0\rangle^{1/4}\int_{[0,T]\times\mathbb{R}} |\overline{v_{\lambda_0}} u_{-\lambda_0-I_{j_1}}u_{-\lambda_0-I_{j_2}}\partial_{x}u_{3\lambda_0+I_{j_3}}| \, dxdt\nonumber\\
\leq &\bigg(\sum_{j_3\approx j_2\approx j_1\geq 0}+\sum_{j_3\approx j_2\gg j_1\geq 0}\bigg)\langle \lambda_0\rangle^{1/4}\int_{[0,T]\times\mathbb{R}} |\overline{v_{\lambda_0}} u_{-\lambda_0-I_{j_1}}u_{-\lambda_0-I_{j_2}}\partial_{x}u_{3\lambda_0+I_{j_3}}| \, dxdt\nonumber\\
:=&I_1(u,v)+I_2(u,v).
\end{align}
For $I_1(u,v)$, from the dispersion modulation decay \eqref{dispersiondecay}, $L^4$ estimate \eqref{lebesgue4a} and Lemma \ref{V2toX}, we have
\begin{align}\label{step2k}
 I_1(u,v) \lesssim & \sum_{j_3\approx j_2\approx j_1\geq 0} \langle \lambda_0\rangle^{1/4}(\langle\lambda_0\rangle+2^{j_3})\|\overline{v_{\lambda_0}}\|_{L^\infty_{x,t}} \| u_{3\lambda_0+I_{j_3}}\|_{L^2_{x,t}} \|u_{-\lambda_0-I_{j_1}}\|_{L^4_{x,t}}\|u_{-\lambda_0-I_{j_2}}\|_{L^4_{x,t}}\nonumber\\
\lesssim & \sum_{j_3\approx j_2\approx j_1\geq 0} \langle \lambda_0\rangle^{1/4}\|v_{\lambda_0}\|_{V^2_A}(\langle\lambda_0\rangle+2^{j_1})
^{-1/2}(\langle\lambda_0\rangle+2^{j_2})^{-1/2}( \langle\lambda_0\rangle+2^{j_3})^{1/4}(2^{j_3})^{1/2}\nonumber\\
&\quad \times T^{\varepsilon/2}(2^{j_1})^{1/4+\varepsilon}(\langle\lambda_0\rangle+2^{j_1})^{-3/8}
(2^{j_2})^{1/4+\varepsilon} (\langle\lambda_0\rangle+2^{j_2})^{-3/8}\|u\|^3_{X^{1/4}_{\infty,A}}\nonumber\\
\lesssim & T^{\varepsilon/2}\sum_{j_3\geq 0} \langle \lambda_0\rangle^{1/4} ( \langle\lambda_0\rangle+2^{j_3})^{-3/2}(2^{j_3})^{1+2\varepsilon}\|v_{\lambda_0}\|_{V^2_A}\|u\|^3_{X^{1/4}_{\infty,A}}\nonumber\\
\lesssim & T^{\varepsilon/2}\bigg(\sum_{j_3\geq 0} (2^{j_3})^{-1/4+2\varepsilon}\bigg) \|v_{\lambda_0}\|_{V^2_A}\|u\|^3_{X^{1/4}_{\infty,A}}\nonumber\\
\lesssim &T^{\varepsilon/2} \|v^{(\lambda_0)}\|_{V^2_A}\|u\|^3_{X^{1/4}_{\infty,A}}.
\end{align}
For $I_2(u,v)$, due to $j_2\gg j_1$, we have $|-\lambda_0-2^{j_2}-\lambda_0-2^{j_1}|\gtrsim (\langle\lambda_0\rangle+2^{j_2})$ and $|-\lambda_0-2^{j_2}+\lambda_0+2^{j_1}|\gtrsim 2^{j_2}$, so we can use bilinear estimate \eqref{bilinear3} to $u_{-\lambda_0-I_{j_1}}$ and $u_{-\lambda_0-I_{j_2}}$. To be specific,
\begin{align}\label{step2l}
 I_2(u,v) \lesssim & \sum_{j_3\approx j_2\gg j_1\geq 0} \langle \lambda_0\rangle^{1/4}(\langle\lambda_0\rangle+2^{j_3})\|\overline{v_{\lambda_0}}\|_{L^\infty_{x,t}} \| u_{3\lambda_0+I_{j_3}}\|_{L^2_{x,t}} \|u_{-\lambda_0-I_{j_1}}u_{-\lambda_0-I_{j_2}}\|_{L^2_{x,t}}\nonumber\\
\lesssim & \sum_{j_3\approx j_2\gg j_1\geq 0} \langle \lambda_0\rangle^{1/4}\|v_{\lambda_0}\|_{V^2_A}(\langle\lambda_0\rangle+2^{j_1})^{-1/2}
(\langle\lambda_0\rangle+2^{j_2})^{-1/2}( \langle\lambda_0\rangle+2^{j_3})^{1/2}\nonumber\\
&\quad \times  (2^{j_3})^{1/2}(\langle\lambda_0\rangle+2^{j_3})^{-1/4} \|u\|_{X^{1/4}_{\infty,A}} T^{\varepsilon/4}(2^{j_2})^{-1/2+\varepsilon}(\langle\lambda_0\rangle+2^{j_2})^{-1/2+\varepsilon}\nonumber\\
&\quad\times(2^{j_1})^{1/2}(\langle\lambda_0\rangle+2^{j_1})^{-1/4}(2^{j_2})^{1/2}
(\langle\lambda_0\rangle+2^{j_2})^{-1/4}\|u\|^2_{X^{1/4}_{\infty,A}}\nonumber\\
\lesssim & T^{\varepsilon/4}\sum_{j_3\gg j_1\geq 0} \langle \lambda_0\rangle^{1/4} (\langle\lambda_0\rangle+2^{j_3})^{-1+\varepsilon}(2^{j_3})^{1/2+\varepsilon}(\langle\lambda_0\rangle+2^{j_1})
^{-1/4}\|v_{\lambda_0}\|_{V^2_A}\|u\|^3_{X^{1/4}_{\infty,A}}\nonumber\\
\lesssim & T^{\varepsilon/4}\bigg(\sum_{j_3\geq 0} (2^{j_3})^{-1/2+2\varepsilon}\cdot j_3\bigg) \|v_{\lambda_0}\|_{V^2_A}\|u\|^3_{X^{1/4}_{\infty,A}}\nonumber\\
\lesssim &T^{\varepsilon/4} \|v^{(\lambda_0)}\|_{V^2_A}\|u\|^3_{X^{1/4}_{\infty,A}}.
\end{align}

If $u_{-\lambda_0-I_{j_2}}$ has the highest dispersion modulation, we still divide the left hand side of \eqref{trilinear4} into two terms as \eqref{step2m}. For $I_1(u,v)$, because of $j_3\approx j_2\approx j_1$, the estimate is exactly same with \eqref{step2k}. For $I_2(u,v)$, we use the bilinear estimate \eqref{bilinear3} to $u_{-\lambda_0-I_{j_1}}$ and $u_{3\lambda_0+I_{j_3}}$. Noticing $|3\lambda_0+2^{j_3}\pm(\lambda_0+2^{j_1})|\gtrsim (\langle\lambda_0\rangle+2^{j_3})$, we have
\begin{align}
 I_2(u,v) \lesssim & \sum_{j_3\approx j_2\gg j_1\geq 0} \langle \lambda_0\rangle^{1/4}(\langle\lambda_0\rangle+2^{j_3})\|\overline{v_{\lambda_0}}\|_{L^\infty_{x,t}} \|u_{-\lambda_0-I_{j_2}}\|_{L^2_{x,t}} \|u_{-\lambda_0-I_{j_1}} u_{3\lambda_0+I_{j_3}}\|_{L^2_{x,t}}\nonumber\\
\lesssim & \sum_{j_3\approx j_2\gg j_1\geq 0} \langle \lambda_0\rangle^{1/4}\|v_{\lambda_0}\|_{V^2_A}(\langle\lambda_0\rangle+2^{j_1})
^{-1/2}(\langle\lambda_0\rangle+2^{j_2})^{-1/2}( \langle\lambda_0\rangle+2^{j_3})^{1/2}\nonumber\\
&\quad \times  (2^{j_2})^{1/2}(\langle\lambda_0\rangle+2^{j_2})^{-1/4} \|u\|_{X^{1/4}_{\infty,A}} T^{\varepsilon/4}(\langle\lambda_0\rangle+2^{j_3})^{-1+2\varepsilon}\nonumber\\
&\quad\times(2^{j_1})^{1/2}(\langle\lambda_0\rangle+2^{j_1})^{-1/4}(2^{j_3})^{1/2}
(\langle\lambda_0\rangle+2^{j_3})^{-1/4}\|u\|^2_{X^{1/4}_{\infty,A}}\nonumber\\
\lesssim & T^{\varepsilon/4}\sum_{j_3\gg j_1\geq 0} \langle \lambda_0\rangle^{1/4} (\langle\lambda_0\rangle+2^{j_3})^{-3/2+2\varepsilon}2^{j_3}(\langle\lambda_0\rangle+2^{j_1})^{-3/4}(2^{j_1})
^{1/2}\|v_{\lambda_0}\|_{V^2_A}\|u\|^3_{X^{1/4}_{\infty,A}}\nonumber\\
\lesssim & T^{\varepsilon/4}\bigg(\sum_{j_3\geq 0} (2^{j_3})^{-1/2+2\varepsilon}\cdot j_3\bigg) \|v_{\lambda_0}\|_{V^2_A}\|u\|^3_{X^{1/4}_{\infty,A}}\nonumber\\
\lesssim &T^{\varepsilon/4} \|v^{(\lambda_0)}\|_{V^2_A}\|u\|^3_{X^{1/4}_{\infty,A}}.
\end{align}

\subsection {$q<\infty$, Proof of \eqref{trilinear3}.}
This subsection $q<\infty$ is similar to the last subsection $q=\infty$, the only difference is to deal with the summation of $\lambda_0$. The frequency constraint condition (FCC) and dispersion modulation constraint condition (DMCC) are same. Thus, we can use the exactly same assortment to $\lambda_0, \cdots, \lambda_3$. Next we take the Case 1 of Step 1 in last subsection for example.

We just denote the left hand side of \eqref{trilinear3} as $\mathscr{L}_{hhhh_-}(u, v)$, and divide it into three parts like the last subsection. For $\mathscr{L}^l_{hhhh_-}(u, v)$, $\lambda_0\approx\lambda_3\approx\lambda_1\approx -\lambda_2$ holds. Thus, From H$\ddot{\rm o}$lder inequality and Strichartz estimate, we have
\begin{align}
 \mathscr{L}^l_{hhhh_{-}}(u, v)&\lesssim\sum_{\lambda_0}\langle \lambda_0\rangle^{5/4} \|\overline{v_{\lambda_0}}\|_{L^4_{x,t}}\|u_{\lambda_0}\|^2_{L^4_{x,t}}\|u_{-\lambda_0}\|_{L^4_{x,t}}\nonumber\\
&\lesssim  T^{1/2}\sum_{\lambda_0}\langle \lambda_0\rangle^{5/4} \|\overline{v_{\lambda_0}}\|_{L^8_tL^4_x}\|u_{\lambda_0}\|^2_{L^8_tL^4_x}\|u_{-\lambda_0}\|_{L^8_tL^4_x}\nonumber\\
&\lesssim  T^{1/2}\sum_{\lambda_0}\langle \lambda_0\rangle^{3/4} \|\overline{v_{\lambda_0}}\|_{V^2_A}\|u_{\lambda_0}\|^2_{U^2_A}\|u_{-\lambda_0}\|_{U^2_A}\nonumber\\
&\lesssim T^{1/2} \|v\|_{Y^0_{q',A}}\|u\|^{3}_{X^{1/4}_{q,A}}. \nonumber
\end{align}

For $\mathscr{L}^m_{hhhh_-}(u, v)$, we still use bilinear estimate \eqref{bilinear3}, Lemma \ref{V2toX}, and H\"older's inequality to obtain that for $0<\varepsilon<1/4q$,
\begin{align*}
\mathscr{L}^m_{hhhh_-}(u, v)&\lesssim\sum_{\lambda_0,j_3\lesssim 1\ll j_1\approx j_2}\langle \lambda_0\rangle^{1/4} \|\overline{v_{\lambda_0}}u_{\lambda_0-I_{j_1}}\|_{L^2_{x,t}}\|u_{-\lambda_0+I_{j_2}}
\partial_{x}u_{\lambda_0-I_{j_3}}\|_{L^2_{x,t}}\nonumber\\
&\lesssim  T^{\varepsilon/2}\sum_{\lambda_0,j_3\lesssim 1\ll j_1\approx j_2}\langle \lambda_0\rangle^{5/4}\langle \lambda_0\rangle^{-1/2+\varepsilon}(2^{j_1})^{-1/2+\varepsilon}\|v_{\lambda_0}\|_{V^2_A} \|u_{\lambda_0-I_{j_1}}\|_{V^2_A}\nonumber\\
&\quad\quad\quad\quad\times \langle \lambda_0\rangle^{-1/2+\varepsilon}(2^{j_2})^{-1/2+\varepsilon}\| u_{-\lambda_0+I_{j_2}}\|_{V^2_A}\|u_{\lambda_0-I_{j_3}}\|_{V^2_A}\nonumber\\
&\lesssim  T^{\varepsilon/2}\sum_{\lambda_0,j_3\lesssim 1\ll j_1\approx j_2}\langle \lambda_0\rangle^{-1/2+2\varepsilon}(2^{j_1})^{\varepsilon-1/q}(2^{j_2})^{\varepsilon-1/q}(2^{j_3})^{1/2-1/q}\nonumber\\
&\quad\quad\quad\quad\quad\quad\times\|v_{\lambda_0}\|_{V^2_A}\|u\|_{X^{1/4}_{q,A}(\lambda_0-I_{j_3})}
\|u\|^2_{X^{1/4}_{q,A}}.
\end{align*}
Making the summation on $j_1,j_2$, then applying H\"older's inequality on $\lambda_0$, and finally summing on $j_3$, we obtain
\begin{align*}
\mathscr{L}^m_{hhhh_-}(u, v)&\lesssim  T^{\varepsilon/2}\sum_{\lambda_0,j_3\leq \log_2\langle\lambda_0\rangle} (2^{j_3})^{4\varepsilon-3/q} \|v_{\lambda_0}\|_{V^2_A}\|u\|_{X^{1/4}_{q,A}(\lambda_0-I_{j_3})}\|u\|^2_{X^{1/4}_{q,A}}\nonumber\\
&\lesssim  T^{\varepsilon/2}\|v\|_{Y^0_{q',A}}\|u\|^{3}_{X^{1/4}_{q,A}}.
\end{align*}

For $\mathscr{L}^h_{hhhh_-}(u, v)$, we just take the case $v_{\lambda_0}$ has the highest dispersion modulation for example and divide $\mathscr{L}^h_{hhhh_-}(u, v)$ into two parts:
\begin{align*}
&\mathscr{L}^{h1}_{hhhh_-}(u, v)+\mathscr{L}^{h2}_{hhhh_-}(u, v)\nonumber\\
:=&\bigg(\sum_{\lambda_0,1\ll j_3\ll j_1\approx j_2}+\sum_{\lambda_0,1\ll j_3\approx j_1\approx j_2}\bigg)\langle \lambda_0\rangle^{1/4}\int_{[0,T]\times\mathbb{R}} |\overline{v_{\lambda_0}} u_{\lambda_0-I_{j_1}}u_{-\lambda_0+I_{j_2}}\partial_{x}u_{\lambda_0-I_{j_3}} |\, dxdt.
\end{align*}
For $\mathscr{L}^{h1}_{hhhh_-}(u, v)$, we have for $0<\varepsilon<1/4q$,
\begin{align*}
\mathscr{L}^{h1}_{hhhh_-}(u, v)&\lesssim \sum_{\lambda_0,1\ll j_3\ll j_1\approx j_2}\langle \lambda_0\rangle^{1/4}\|\overline{v_{\lambda_0}}\|_{L^2_tL^{\infty}_x} \|u_{\lambda_0-I_{j_1}}\|_{L^\infty_tL^{2}_x}\| u_{-\lambda_0+I_{j_2}}\partial_{x}u_{\lambda_0-I_{j_3}}\|_{L^2_{x,t}}\nonumber\\
&\lesssim \sum_{\lambda_0,1\ll j_3\ll j_1\approx j_2} \langle \lambda_0\rangle^{5/4}\langle \lambda_0\rangle^{-1/2}(2^{j_1})^{-1/2}(2^{j_3})^{-1/2}\|v_{\lambda_0}\|_{V^2_A} \|u_{\lambda_0-I_{j_1}}\|_{V^2_A}\nonumber\\
&\quad\quad \times T^{\varepsilon/4}\langle \lambda_0\rangle^{-1/2+\varepsilon}(2^{j_2})^{-1/2+\varepsilon}\| u_{-\lambda_0+I_{j_2}}\|_{V^2_A}\|u_{\lambda_0-I_{j_3}}\|_{V^2_A}\nonumber\\
&\lesssim  T^{\varepsilon/4}\sum_{\lambda_0,1\ll j_3\ll j_1\approx j_2}\langle \lambda_0\rangle^{-1/2+\varepsilon}(2^{j_2})^{\varepsilon-2/q}(2^{j_3})^{-1/q}
\|v_{\lambda_0}\|_{V^2_A}\|u\|_{X^{1/4}_{q,A}(\lambda_0-I_{j_3})}
\|u\|^2_{X^{1/4}_{q,A}}\nonumber\\
&\lesssim  T^{\varepsilon/4}\sum_{\lambda_0,1\ll j_3\leq\log_2\langle\lambda_0\rangle}(2^{j_3})^{-1/2+2\varepsilon-3/q}
\|v_{\lambda_0}\|_{V^2_A}\|u\|_{X^{1/4}_{q,A}(\lambda_0-I_{j_3})}
\|u\|^2_{X^{1/4}_{q,A}}\nonumber\\
&\lesssim  T^{\varepsilon/4}\|v\|_{Y^0_{q',A}}\|u\|^{3}_{X^{1/4}_{q,A}}.
\end{align*}
For $\mathscr{L}^{h2}_{hhhh_-}(u, v)$,  we know that for $0<\varepsilon<1/4q$,
\begin{align*}
\mathscr{L}^{h2}_{hhhh_-}(u, v)&\lesssim \sum_{\lambda_0,1\ll j_3\approx j_1\approx j_2}\langle \lambda_0\rangle^{5/4}\|\overline{v_{\lambda_0}}\|_{L^2_tL^{\infty}_x} \|u_{\lambda_0-I_{j_1}}\|_{L^\infty_tL^{2}_x}\| u_{-\lambda_0+I_{j_2}}\|_{L^4_{x,t}}\|u_{\lambda_0-I_{j_3}}\|_{L^4_{x,t}}\nonumber\\
&\lesssim \sum_{\lambda_0,1\ll j_3\approx j_1\approx j_2}\langle \lambda_0\rangle^{5/4}\langle \lambda_0\rangle^{-1/2}(2^{j_1})^{-1/2}(2^{j_3})^{-1/2}\|v_{\lambda_0}\|_{V^2_A} \|u_{\lambda_0-I_{j_1}}\|_{V^2_A}\nonumber\\
&\quad\quad \times T^{\varepsilon/2}\langle \lambda_0\rangle^{-3/4}(2^{j_2})^{1/4-1/q+\varepsilon}(2^{j_3})^{1/4-1/q+\varepsilon}\| u\|^2_{X^{1/4}_{q,A}}\nonumber\\
&\lesssim T^{\varepsilon/2}\sum_{\lambda_0,j_1}(2^{j_1})^{2\varepsilon-3/q}
\|v_{\lambda_0}\|_{V^2_A}\|u\|_{X^{1/4}_{q,A}(\lambda_0-I_{j_1})}\|u\|^2_{X^{1/4}_{q,A}}\nonumber\\
&\lesssim T^{\varepsilon/2}\|v\|_{Y^0_{q',A}}\|u\|^{3}_{X^{1/4}_{q,A}}.
\end{align*}
Where the last inequality is by applying H\"older's inequality on $\lambda_0$. For other cases, we can take the similar calculation to get the desired estimates, thus we omit it.

\section{Ill-posedness result}

In this section we study the Cauchy problem of the defocusing mKdV equation ( the focusing case can also be treated by our method):
\begin{align}\label{illposedPDE}
   u_{t} + u_{xxx} - (u^3)_x = 0,\quad
u(0,x)= \delta u_0.
\end{align}
We have the ill-posedness result as follows.
\begin{thm}\label{illtheorem} Let $s<1/4$, $2\leq q\leq \infty$, $0<\delta\ll 1$. Then for the mKdV equation \eqref{illposedPDE}, the solution map $\delta u_0\rightarrow u(\delta,t)$ in  $M^s_{2,q}$ is not $C^3$ continuous at the origin.
\end{thm}

{\bf Proof.} From \eqref{ImKdV} we can define the solution map as follows:
\begin{align}
\mathcal{T}: \ \ \delta u_0\rightarrow u(\delta,t)=e^{-t\partial_x^3}\delta u_0 +\int^t_0 e^{-(t-\tau)\partial_x^3}{(u^3)_x}(\tau)\, d\tau. \label{integral}
\end{align}
By straightforward calculations, we get
 \begin{align}
   u(\delta, t)|_{\delta=0} =0;&\ \ \ \ u_1:=  \frac{\partial u}{\partial\delta}\bigg|_{\delta=0}= e^{-t\partial_x^3}u_0;\  \ \  \  u_2:= \frac{\partial^2 u}{\partial\delta^2}\bigg|_{\delta=0}=0;\label{u1}\\
   u_3:= & \frac{\partial^3 u}{\partial\delta^3}\bigg|_{\delta=0}=6 \int^t_0 e^{-(t-\tau)\partial_x^3}\partial_x(e^{-\tau\partial_x^3}u_0)^3\, d\tau. \label{u3}
 \end{align}

It is well known that if the map $\delta u_0\rightarrow u(\delta)$ is of class $C^3$  at the origin, the necessary condition is
\begin{align}\label{illposed0}
\sup_{t\in[0,T]}\|u_3\|_{M^s_{2,q}}\leq C\|u_0\|_{M^s_{2,q}}^3.
\end{align}

We choose a suitable $u_0\in M^s_{2,q}$, $s<1/4$ defined by  $$\widehat{u_0}(\xi)=N^{-s+1/4}\big(\chi_{[N,N+\frac{1}{\sqrt{N}}]}(\xi)+\chi_{[-N-\frac{1}{\sqrt{N}},-N]}(\xi)\big).$$
Note that $\|u_0\|_{M^s_{2,q}}\sim1$.

We estimate the Fourier transform of $u_3$ in \eqref{u3} as follows
\begin{align}\label{illposed1}
  \widehat{u_3}(\xi) \simeq & \int^t_0 e^{{\rm i}(t-\tau)\xi^3}({\rm i}\xi)\int_{\mathbb{R}^2}e^{{\rm i}\tau(\xi-\xi_1-\xi_2)^3}\widehat{u_0}(\xi-\xi_1-\xi_2)e^{{\rm i}\tau\xi_1^3}\widehat{u_0}(\xi_1)e^{{\rm i}\tau\xi_2^3}\widehat{u_0}(\xi_2)d\xi_1d\xi_2d\tau\nonumber\\
  \simeq&e^{{\rm i}t\xi^3}({\rm i}\xi)\int_{\mathbb{R}^2}\int^t_0 e^{{\rm i}\tau\Phi(\xi,\xi_1,\xi_2)}d\tau\, \widehat{u_0}(\xi-\xi_1-\xi_2)\widehat{u_0}(\xi_1)\widehat{u_0}(\xi_2)d\xi_1d\xi_2\nonumber\\
  \simeq&e^{{\rm i}t\xi^3}({\rm i}\xi)\int_{\mathbb{R}^2}\frac{e^{{\rm i}t\Phi(\xi,\xi_1,\xi_2)}-1}{{\rm i}\Phi(\xi,\xi_1,\xi_2)}\widehat{u_0}(\xi-\xi_1-\xi_2)\widehat{u_0}(\xi_1)\widehat{u_0}(\xi_2)d\xi_1d\xi_2,
\end{align}
where $\Phi(\xi,\xi_1,\xi_2)=-3(\xi-\xi_1)(\xi-\xi_2)(\xi_1+\xi_2)$. Noticing that for $\xi-\xi_1-\xi_2$, $\xi_1$, and $\xi_2$,   if one or two items of them locate in $[N,N+1/\sqrt{N}]$, we have $|\Phi(\xi,\xi_1,\xi_2)|\lesssim 1$; if all three items locate in $[N,N+1/\sqrt{N}]$ (or $[-N-1/\sqrt{N},-N]$), we have $|\Phi(\xi,\xi_1,\xi_2)|\sim N^3$, then \eqref{illposed1} shall be much smaller. Therefore,
\begin{align}
  \widehat{u_3}(\xi)\simeq&N^{-3s+3/4}e^{{\rm i}t\xi^3}({\rm i}\xi)\int_{\mathbb{R}^2}\frac{e^{{\rm i}t\Phi(\xi,\xi_1,\xi_2)}-1}{{\rm i}\Phi(\xi,\xi_1,\xi_2)}\nonumber\\
&\quad\quad\times\chi_{[N,N+\frac{1}{\sqrt{N}}]}(\xi-\xi_1-\xi_2)\chi_{[N,N+\frac{1}{\sqrt{N}}]}(\xi_1)
\chi_{[-N-\frac{1}{\sqrt{N}},-N]}(\xi_2)d\xi_1d\xi_2\nonumber\\
\simeq&N^{-3s+3/4}e^{{\rm i}t\xi^3}\xi\int_{N}^{N+\frac{1}{\sqrt{N}}}\int_{-N-\frac{1}{\sqrt{N}}}^{-N}    te^{{\rm i}t\theta}\cdot\chi_{[N,N+\frac{1}{\sqrt{N}}]}(\xi-\xi_1-\xi_2)d\xi_1d\xi_2,\nonumber
\end{align}
where $\theta\in [0,\Phi]$ or $[\Phi,0]$, $|\Phi(\xi,\xi_1,\xi_2)|=O(1)$, $\xi\in[N-1/\sqrt{N},N+2/\sqrt{N}]$. Thus there exists  a small and fixed constant $t$ such that
\begin{align*}
\|u_3\|_{M^s_{2,q}}\geq C N^{(-2s+1/2)}\ \ \  (2\leq q\leq\infty).
\end{align*}
We find that  \eqref{illposed0} leads to
\begin{align*}
  -2s+1/2\leq0\ \ {\rm  i.e.}\ \ s\geq 1/4.
\end{align*}
Now we complete the proof. $\hfill\Box$

\noindent{\bf Acknowledgments.} The author is indebted to Prof. Baoxiang Wang for very helpful discussions, encouragements and supports. This work is supported in part by the National Science Foundation of China, grants 11271023 and 11571254.

\end{document}